\documentclass[12pt] {article}
\usepackage{amsmath,amsthm}
\usepackage{amssymb,latexsym}
\usepackage{graphicx}
\usepackage{amssymb}
\usepackage{amsfonts}
\usepackage{epic,eepic}
\usepackage{amscd}

\usepackage{amscd}
\title{Plurisubharmonic functions on the octonionic plane and $Spin(9)$-invariant valuations on convex sets.}
\date{}
\author{ Semyon Alesker \footnote{Partially supported by ISF grant 1369/04.}
\\  { \normalsize Department of Mathematics, Tel Aviv University, Ramat Aviv}
\\  { \normalsize 69978 Tel Aviv,
Israel }
\\ {\normalsize e-mail: semyon@post.tau.ac.il}}

\def\eps{\varepsilon}
\def\alp{\alpha}

\def\Ome{\Omega}
\def\lam{\lambda}

\def\to{\rightarrow}

\def\qed { Q.E.D. }
\def\pt{\partial}

\def\RR{\mathbb{R}}
\def\CC{\mathbb{C}}

\def\NN{\mathbb{N}}

\def\HH{\mathbb{H}}

\def\PP{\mathbb{P}}
\def\OO{\mathbb{O}}

\def\dfq{\frac{\partial ^2 f}{\partial\bar q_i \partial q_j}}
\def\duq{\frac{\partial ^2 u}{\partial\bar q_i \partial q_j}}

\def\ddbq{\frac{\partial ^2 b}{\partial \bar q_i \partial q_j}}

\swapnumbers
\newtheorem{theorem}{Theorem}[subsection]
\newtheorem{corollary}[theorem]{Corollary}

\newtheorem{lemma}[theorem]{Lemma}
\newtheorem{proposition}[theorem]{Proposition}
\newtheorem{claim}[theorem]{Claim}
\theoremstyle{definition}
\newtheorem{example}[theorem]{Example}
\newtheorem{definition}[theorem]{Definition}
\newtheorem{remark}[theorem]{Remark}

\def\ca{{\cal A}}  
\def\cd{{\cal D}}  
 \def\ch{{\cal H}} 
 \def\ck{{\cal K}} 
  
\def\cp{{\cal P}}

\def\inj{\hookrightarrow}

\newcommand \supp{\operatorname{supp} \,}

\begin{document}
\maketitle \setcounter{section}{-1}
\begin{abstract}
A new class of plurisubharmonic functions on the octonionic plane
$\OO^2\simeq\RR^{16}$ is introduced. An octonionic version of
theorems of A.D. Aleksandrov \cite{aleksandrov-58} and
Chern-Levine-Nirenberg \cite{chern-levine-nirenberg}, and B{\l}ocki
\cite{blocki} are proved. These results are used to construct new
examples of continuous translation invariant valuations on convex
subsets of $\OO^2\simeq \RR^{16}$. In particular a new example of
$Spin(9)$-invariant valuation on $\RR^{16}$ is given.
\end{abstract}
\tableofcontents
\section{Introduction}\label{S:intro}
\subsection{An overview}\label{Ss:overview}
Let $\OO$ denote the division algebra of octonions (=Cayley
numbers). The goal of this article is to introduce and study a new
class of plurisubharmonic functions on the octonionic plane
$\OO^2\simeq\RR^{16}$, and to give some applications to the theory
of valuations on convex sets. In particular octonionic versions of
theorems of A.D. Aleksandrov \cite{aleksandrov-58} and
Chern-Levine-Nirenberg \cite{chern-levine-nirenberg}, and B{\l}ocki
\cite{blocki} are proved. These results are used to construct new
examples of continuous translation invariant valuations on convex
subsets of $\OO^2\simeq \RR^{16}$. In particular a new example of a
$Spin(9)$-invariant valuation on $\RR^{16}$ is given (notice that
$Spin(9)$ is one of the three exceptional examples of compact
connected groups acting transitively on a sphere - see discussion
below in the introduction).

The theories of convex functions on $\RR^n$ and plurisubharmonic
functions on complex manifolds are classical and well studied: see
e.g. the book \cite{hormander} for the introduction to these
subjects; the book by Lelong \cite{lelong} is a classical
introduction to the theory of plurisubharmonic functions of complex
variables.


More recently a class of plurisubharmonic functions of quaternionic
variables on the flat quaternionic space $\HH^n$ has been introduced
independently and at the same time by the author and G. Henkin. It
was investigated further and applied by the author
\cite{alesker-busm-03}, \cite{alesker-jga-03},\cite{alesker-adv-05}
and G. Henkin \cite{henkin}.
Then part of this theory has been generalized to more general
context of (not necessarily flat) hypercomplex manifolds by M.
Verbitsky and the author \cite{alesker-verbitsky}. We refer also
to \cite{alesker-alg-anal} for a survey of these results. Very
recently, other classes of plurisubharmonic functions have been
introduced in the context of calibrated geometries
\cite{harvey-lawson-psh}.

Let us describe our main results in greater details. The algebra
$\OO$ of octonions is a non-associative non-commutative division
algebra over the reals $\RR$ of dimension 8. $\OO$ has a basis over
$\RR$: $e_0,e_1,\dots,e_7$ where $e_0=1$ is the identity element,
and for $i>0$ the $e_i$'s are anti-commuting elements satisfying
$e_i^2=-1$. Any octonion $q\in \OO$ can be uniquely written in the
form
$$q=\sum_{i=0}^7x_ie_i \mbox{ with } x_i\in \RR.$$
One introduces an {\itshape octonionic conjugation} by $q\mapsto
\bar q:=x_0-\sum_{i=1}^7x_ie_i$. This conjugation is an
anti-involution on $\OO$. These and some other elementary properties
of $\OO$ are reviewed in more details in Section \ref{S:oct-prop}.
For further properties we refer to the survey article \cite{baez}
and Chapters 6,14 of the book \cite{harvey}.

Let $F\colon \OO\to \OO$ be a smooth function. We introduce two
{\itshape Dirac operators} as follows
\begin{eqnarray}\label{Dirac-bar}
\frac{\pt}{\pt\bar q}F:=\sum_{i=0}^7e_i\frac{\pt F}{\pt
x_i},\\\label{Dirac} \frac{\pt}{\pt q}F:=\overline{\frac{\pt \bar
F}{\pt \bar q}}=\sum_{i=0}^7\frac{\pt F}{\pt x_i}\bar e_i.
\end{eqnarray}
\begin{remark}
These operators have been introduced in analogy to the complex and
quaternionic cases where they are well known (see
\cite{alesker-busm-03} and references therein).
\end{remark}

Similarly if $F\colon \OO^m\to \OO$ is a smooth function of $m$
octonionic variables then one can define operators $\frac{\pt}{\pt
q_i},\, \frac{\pt}{\pt\bar q_i}$ for $i=1,\dots,m$.  It is easy to
see that if the function $F$ is real valued then for any
$i,j=1,\dots,m$
$$\frac{\pt}{\pt q_i}\left(\frac{\pt}{\pt\bar q_j}F\right)=\frac{\pt}{\pt\bar
q_j}\left(\frac{\pt}{\pt q_i}F\right).$$ This expression will be
denoted either by $\frac{\pt^2F}{\pt q_i\pt\bar q_j}$ or by
$\frac{\pt^2F}{\pt \bar q_j\pt q_i}$. The matrix with octonionic
entries $\left(\frac{\pt^2F}{\pt\bar q_i\pt q_j}\right)_{i,j=1}^m$
will be called the {\itshape octonionic Hessian} of a real valued
function $F$.

In fact the octonionic Hessian of a real valued function is an
octonionic hermitian matrix. A matrix $A=(a_{ij})$ is called
octonionic hermitian if
$$a_{ij}=\overline{a_{ji}} \mbox{ for any } i,j.$$
Let us discuss now the case of $m=2$ octonionic variables (the case
$m=1$ is trivial, the case $m\geq 3$ seems to be quite different
from the case $m=2$ and it has not been studied). Let $\Ome\subset
\OO^2$ be an open subset.
\begin{definition}\label{D:psh-intro}
A function $f\colon \Ome\to \RR\cup\{-\infty\}$ is called {\itshape
octonionic plurisubharmonic} if $f$ is upper semi-continuous and its
restrictions to any affine octonionic line is subharmonic.
\end{definition}
Affine octonionic lines are discussed in Section \ref{S:radon}
below; relevant notions used in Definition \ref{D:psh-intro} are
recalled in Section \ref{Ss:psh-func}.

It is shown in Proposition \ref{P:smooth-psh} that a twice
continuously differentiable function $f$ is octonionic
plurisubharmonic if and only if the octonionic Hessian
$\left(\frac{\pt^2f}{\pt\bar q_i\pt q_j}\right)_{i,j=1}^2$ is
non-negative definite octonionic hermitian matrix pointwise in the
sense of Definition \ref{D:positive-def} below. Notice that any
convex function is octonionic plurisubharmonic. Any octonionic
plurisubharmonic function is subharmonic (Proposition
\ref{P:psh-sh}). The sum and the maximum of finitely many
octonionic plurisubharmonic functions are octonionic
plurisubharmonic (Proposition \ref{P:sums-maximums}). The class of
octonionic plurisubharmonic functions is invariant under
translations and linear transformations from the group
$SL_2(\OO)\simeq Spin(9,1)$ (which is discussed in Section
\ref{Ss:sl2}). We denote by $P(\Ome)$ the class of all octonionic
plurisubharmonic functions in $\Ome$.

On the class of octonionic hermitian $(2\times 2)$-matrices there
exists a determinant function $\det$ with various nice properties
completely analogous to properties of the usual determinant of
real symmetric and complex hermitian matrices (and also of the
Moore determinant of quaternionic hermitian matrices). We refer to
Section \ref{S:la} for the definition and the main properties of
it. This determinant plays an important role in our constructions.
Note also that for octonionic hermitian matrices of size at least
4 no nice notion of determinant is known, while for matrices of
size 3 it does exist (see e.g. Section 3.4 in \cite{baez}).

The first main result can be stated as follows (see Proposition
\ref{P:cln0} and Theorem \ref{T:cln-main}).
\begin{theorem}\label{T:cln-intro}
For any $f\in P(\Ome)\cap C(\Ome)$ there exists a non-negative
measure in $\Ome$ denoted by $\det\left(\dfq\right)$ which is
uniquely characterized by the following two properties:

(i) this measure has the obvious meaning if $f\in C^2(\Ome)$;

(ii) if a sequence $\{f_n\}\subset P(\Ome)\cap C(\Ome)$ converges
uniformly on compact subsets to the function $f$ then
$$\det\left(\frac{\pt^2 f_n}{\pt\bar q_i\pt q_j}\right)\to
\det\left(\dfq\right)$$ weakly in sense of measures.
\end{theorem}
\begin{remark}
A real version of this result for the usual Hessian of convex
functions was proved by A.D. Aleksandrov \cite{aleksandrov-58}. A
complex version for the complex Hessian of complex plurisubharmonic
functions was proved by Chern-Levine-Nirenberg
\cite{chern-levine-nirenberg}. A quaternionic version for
quaternionic Hessian of quaternionic plurisubharmonic functions was
proved by the author \cite{alesker-busm-03} on the flat space
$\HH^n$ and by M. Verbitsky and the author \cite{alesker-verbitsky}
on more general hypercomplex manifolds.
\end{remark}
The second main result is the following octonionic version of a
result by B{\l}ocki \cite{blocki} for complex plurisubharmonic
functions (a quaternionic version was proved by the author in
\cite{alesker-adv-05}).
\begin{theorem}\label{T:blocki-intro}
For any $u,v\in P(\Ome)\cap C(\Ome)$ such that $\min\{u,v\}\in
P(\Ome)$ one has
$$\det(\pt^2(\min\{u,v\}))=\det(\pt^2u)+\det(\pt^2v)-\det(\pt^2(\max\{u,v\})).$$
\end{theorem}
Actually this result is an easy consequence of a more precise
Theorem \ref{T:blocki1} (which is due to B{\l}ocki \cite{blocki} in
the complex case).

We refer to Section \ref{S:psh} for other results on octonionic
plurisubharmonic functions. Now we are going to describe
applications of the above results to the theory of valuations on
convex sets. First let us remind basic notions of this theory
referring for further information to the surveys by McMullen
\cite{mcmullen-survey} and McMullen and Schneider
\cite{mcmullen-schneider}. Let $V$ be a finite dimensional real
vector space. Let ${\cal K}(V)$ denote the class of all non-empty
convex compact subsets of $V$.

\begin{definition}
(i) A function $\phi :{\cal K}(V) \to \CC$ is called a {\itshape
valuation} if for any $K_1, \, K_2 \in {\cal K}(V)$ such that their
union is also convex one has
$$\phi(K_1 \cup K_2)= \phi(K_1) +\phi(K_2) -\phi(K_1 \cap K_2).$$

(ii) A valuation $\phi$ is called {\itshape continuous} if it is
continuous with respect the Hausdorff metric on ${\cal K}(V)$.
\end{definition}

Recall that the Hausdorff metric $d_H$ on ${\cal K}(V)$ depends on
a choice of a Euclidean metric on $V$ and it is defined as
follows: $d_H(A,B):=\inf\{ \eps >0|A\subset (B)_\eps \mbox{ and }
B\subset (A)_\eps\}$ where $(U)_\eps$ denotes the
$\eps$-neighborhood of a set $U$. Equipped with the Hausdorff
metric, ${\cal K} (V)$ becomes a locally compact space, and the
induced topology on $\ck(V)$ is independent of a choice of the
Euclidean metric on $V$. Let us denote by $Val(V)$ the space of
translation invariant continuous valuations on $V$.

\begin{example}\label{ex-valuations}
(1) A Lebesgue measure $vol$ on $V$ belongs to $Val(V)$.

(2) The Euler characteristic $\chi$ belongs to $Val(V)$. Recall that
$\chi (K)=1$ for any $K\in \ck(V)$.

(3) Denote $m:=\dim V$. Fix $k=1,\dots,m$. Fix $A_1,\dots,A_{m-k}\in
\ck(V)$. Then the mixed volume $$K\mapsto
V(K[k],A_1,\dots,A_{m-k})$$ belongs to $Val(V)$ (here $K[k]$ means
that a set $K$ is taken $k$ times). For the notion of mixed volume
and its properties see e.g. the book \cite{schneider-book}.
\end{example}
It was conjectured by P. McMullen \cite{mcmullen-80} and proved by
the author \cite{alesker-gafa-01} that the linear combinations of
the mixed volumes as in Example \ref{ex-valuations} (3) above are
dense in $Val(V)$ in the topology of uniform convergence on compact
subsets of $\ck(V)$. Nevertheless there are other than mixed volumes
non-trivial constructions of translation invariant continuous
valuations. One of such constructions will be described in Theorem
\ref{T:valuat-intro} below, which provides in particular a new
example of a continuous $Spin(9)$-invariant valuation on $\RR^{16}$.


To explain why such examples are interesting let us digress to a
more general context of valuations invariant under a group. Assume
that $V$ is a Euclidean space. Let $G$ be a compact subgroup of the
orthogonal group. Let us denote by $Val^G(V)$ the space of
$G$-invariant continuous translation invariant valuations. It is
known \cite{alesker-adv-00} that the space $Val^G(V)$ is finite
dimensional if and only if $G$ acts transitively on the unit sphere
in $V$. In such a case one can try to classify explicitly this
space.

In topology there exists a complete classification of compact
connected Lie groups acting transitively on spheres
\cite{montgomery-samelson}, \cite{borel1}, \cite{borel2}. It is
shown that there exist

$\bullet$ 6 infinite series: $SO(n), U(n), SU(n), Sp(n),Sp(n)\cdot
Sp(1), Sp(n)\cdot U(1)$;

$\bullet$ 3 exceptions: $G_2, Spin(7), Spin(9)$.

If $G$ is either the full orthogonal or special orthogonal group
the corresponding classification of $G$-invariant valuations is
well known and this is the famous result by Hadwiger
\cite{hadwiger-book}. The classification for the unitary group
$U(n)$ acting on $\CC^n\simeq \RR^{2n}$ was obtained by the author
\cite{alesker-jdg-03} (see Fu's article \cite{fu-jdg-06} for
further information on the algebra structure of
$Val^{U(n)}(\CC^n)$). The case of $G=SU(2)$ has been classified by
the author \cite{alesker-su2} (see Bernig's article
\cite{bernig-su2} for the algebra structure of
$Val^{SU(2)}(\CC^2)$).

Except of these groups $O(n), SO(n), U(n),$ and $SU(2)$ no
classification of valuations has been obtained so far. To obtain
such a classification is an interesting question which does not
seem to have an easy solution. In \cite{alesker-adv-05} some
non-trivial examples of valuations invariant under the
quaternionic groups have been constructed. The motivation of this
article is to say something on the exceptional group $Spin(9)$
which acts transitively on the sphere $S^{15}\subset \RR^{16}$.
This group is discussed in Section \ref{Ss:projline}, see in
particular Remark \ref{R:spin9}. It will be convenient to identify
$\RR^{16}$ with the octonionic plane $\OO^2$.


The next main result of this article is a new construction of
continuous valuations on $\OO^2\simeq \RR^{16}$ based on
octonionic plurisubharmonic functions. For  a convex set $K\in
\ck(\OO^2)$ we denote its supporting functional by $h_K\colon
\OO^2\to \RR$ (the definition is recalled in Section
\ref{Ss:new-examples}).
\begin{theorem}\label{T:valuat-intro}
Fix a continuous compactly supported function $\psi$ on $\OO^2$.
Then
$$K\mapsto \int_{\OO^2}\det\left(\frac{\pt^2h_K}{\pt\bar q_i\pt q_j}\right)\cdot
\psi dq$$ is a translation invariant continuous valuation on
$\ck(\OO^2)$.
\end{theorem}
The continuity is a consequence of Theorem \ref{T:cln-intro}, and
the valuation property is a consequence of Theorem
\ref{T:blocki-intro}. As an immediate corollary we get that
$$P_\OO(K):=\int_D\det\left(\frac{\pt^2h_K}{\pt\bar q_i\pt q_j}\right) dq,$$
where $D$ is the unit centered ball in $\OO^2$, is a continuous
translation invariant $Spin(9)$-invariant valuation. We call $P_\OO$
the {\itshape octonionic pseudo-volume}. Other examples of
$Spin(9)$-invariant valuations of different nature coming from the
convex and integral geometry are described in Section
\ref{Ss:obvious-examples}.

It should be noted that a complex (and in fact the original) version
of the pseudo-volume using the complex Hessian was first considered
in the context of convexity (though not of valuations) by
Kazarnovski\u\i \cite{kazarnovskii-81}, \cite{kazarnovskii-84}. The
quaternionic version of the pseudo-volume using the quaternionic
Hessian was constructed by the author in \cite{alesker-adv-05}. As a
side remark notice that the real version of the pseudo-volume is
proportional to the usual volume.

{\bf Acknowledgements.} I thank J. Bernstein, M. Borovoi, and M.
Sodin for useful conversations. Part of this work was done during my
stay at the Institute for Advanced Study at Jerusalem; I thank this
institute for the hospitality.

\subsection{Organization on the article}\label{Ss:organization}
The article is organized as follows. In Section \ref{S:oct-prop} we
collect various definitions and facts related to octonions including
some linear algebra over octonions. Probably no result of this
section is new (with the only possible exception of Proposition
\ref{P:ohs1}). A reader familiar with octonions can skip this
section and consult it only whenever necessary.

In Section \ref{S:radon} we prove that the Radon transform over the
set of affine octonionic lines in $\OO^2$ is injective. We prove it
by constructing an explicit inversion formula. Most probably this
result is not new. We need it for some technical reasons (proof of
Lemma \ref{L:psh1}) and present a proof due to the lack of a
reference.

Sections \ref{S:psh} and \ref{S:valuat} are the main ones. In
Section \ref{S:psh} we introduce the class of plurisubharmonic
functions on the octonionic plane $\OO^2$ and establish our main
results on this class.

In Section \ref{S:valuat} we discuss in  detail applications of the
above technique to valuations on convex sets in $\OO^2\simeq
\RR^{16}$ and describe some other $Spin(9)$-invariant valuations..

\subsection{Notation list.}\label{Ss:notation}
The following notation will be used often throughout the article.

$\bullet$ $\HH$, $\OO$ - algebras of quaternions and  octonions
respectively;

$\bullet$ $\left(\duq\right)$ or $\left(\pt^2 u\right)$ - the
octonionic Hessian of a real valued function $u$.

$\bullet$ $C(\Ome)$ the space of continuous functions in a domain
$\Ome$.

$\bullet$ $C^\infty(\Ome)$ (resp. $C^{\infty}_0(\Ome)$) - the space
of complex valued infinitely smooth functions (resp. with compact
support) in $\Ome$.

$\bullet$ $C^{-\infty}(\Ome)$ - the space of complex valued
generalized functions in $\Ome$ (by definition, this space is the
topological dual to the space of infinitely smooth densities with
compact support in $\Ome$).

$\bullet$ $L^1_{loc}(\Ome)$ - the space of locally integrable
functions in $\Ome$.

$\bullet$ $P(\Ome)$ - the space of octonionic plurisubharmonic
functions in $\Ome\subset \OO^2$.

$\bullet$ $\OO\PP^1$ - the octonionic projective line (see Section
\ref{Ss:projline}).

$\bullet$ $\ck(V)$ - the family of non-empty convex compact subsets
of a vector space $V$.

$\bullet$ $Val(V)$ - the space of translation invariant continuous
valuations on convex compact subsets of $V$.

$\bullet$ $Val^G(V)$ - the space of translation invariant continuous
valuations on $V$ which are invariant under a group $G$.

\section{Basic properties of the octonions.}\label{S:oct-prop}
In this section we collect various facts on octonions. Probably no
result in this section is new with the only possible exception of
Proposition \ref{P:ohs1}. Whenever possible we give references.
Otherwise we present a proof. The reader is advised to consult the
survey article \cite{baez} and Chapters 6,14 of the book
\cite{harvey}.

\subsection{Some octonionic algebra.}\label{Ss:oalgebra}
The octonions $\OO$ form an 8-dimensional algebra over the reals
$\RR$ which is neither associative nor commutative. The product can
be described as follows. $\OO$ has a basis $e_0,e_1,\dots,e_7$ over
$\RR$ where $e_0=1$ is the unit, and the product can be given by the
multiplication table:

\begin{center}
$$\begin{array}{c|c|c|c|c|c|c|c|}
   & e_1& e_2& e_3& e_4& e_5& e_6& e_7\\\hline
e_1&-1  &e_4 &e_7 &-e_2&e_6 &-e_5&-e_3\\\hline e_2&-e_4&-1
&e_5&e_1 &-e_3 &e_7&-e_6\\\hline
e_3&-e_7&-e_5&-1&e_6&e_2&-e_4&e_1\\\hline
e_4&e_2&-e_1&-e_6&-1&e_7&e_3&-e_5\\\hline
e_5&-e_6&e_3&-e_2&-e_7&-1&e_1&e_4\\\hline
e_6&e_5&-e_7&e_4&-e_3&-e_1&-1&e_2\\\hline
e_7&e_3&e_6&-e_1&e_5&-e_4&-e_2&-1\\\hline
\end{array}$$
\end{center}

There is also another easier way to remember the product using so
called {\itshape Fano plane} (see Figure 1 below). In Figure 1 each
pair of distinct points lies in a unique line (the circle is also
considered to be a line). Each line contains exactly three points,
and these points are cyclically oriented. If $e_i,e_j,e_k$ are
cyclically oriented in this way then
$$e_ie_j=-e_je_i=e_k.$$
We have to add two more rules:

$\bullet$ $e_0=1$ is the identity element;

$\bullet$ $e_i^2=-1$ for $i>0$.

All these rules define uniquely the algebra structure of $\OO$. The
center of $\OO$ is equal to $\RR$.



\begin{figure}[h]
\setlength{\unitlength}{0.00087489in}
\begingroup\makeatletter\ifx\SetFigFont\undefined%
\gdef\SetFigFont#1#2#3#4#5{%
  \reset@font\fontsize{#1}{#2pt}%
  \fontfamily{#3}\fontseries{#4}\fontshape{#5}%
  \selectfont}%
\fi\endgroup%
{\renewcommand{\dashlinestretch}{30}
\begin{picture}(9591,4295)(-1541,1000)
\path(9409,957)(9454,1137)(9499,957)
\path(2145,3620)(2101,3441)(2057,3620)
\path(2145,2405)(2101,2226)(2057,2405)
\path(1061,1673)(881,1710)(1062,1764)
\path(3326,1658)(3146,1710)(3327,1749)
\path(508,2415)(633,2549)(586,2370)
\path(1533,4196)(1672,4320)(1651,4151)
\path(2584,4349)(2620,4170)(2468,4302)
\path(3624,2564)(3666,2385)(3538,2517)
\path(938,2580)(987,2766)(1065,2621)
\path(1599,3070)(1478,3228)(1661,3193)
\path(2589,2530)(2468,2648)(2641,2593)
\path(1515,2564)(1688,2620)(1567,2490)
\path(2430,3104)(2603,3165)(2482,3030)
\path(2416,3957)(2559,3860)(2378,3875)
\path(2918,2012)(2758,1930)(2849,2070)
\put(2100,2837){\ellipse{2250}{2250}}
\put(4066,1710){\whiten\ellipse{346}{346}}
\put(2100,1710){\whiten\ellipse{346}{346}}
\put(2100,2850){\whiten\ellipse{346}{346}}
\put(3056,3422){\whiten\ellipse{346}{346}}
\put(1146,3422){\whiten\ellipse{346}{346}}
\put(141,1710){\whiten\ellipse{346}{346}}
\put(2100,5047){\whiten\ellipse{346}{346}}
\path(4066,1710)(2100,5047)
\path(2100,5047)(2100,1710) 
\path(141,1710)(4066,1710)
\path(141,1710)(2100,5047)
\path(4066,1710)(1146,3422)
\path(141,1710)(3056,3422)
\put(98,1630){\makebox(0,0)[lb]{\smash{{\SetFigFont{12}{14.4}{\rmdefault}{\mddefault}{\updefault}$\scriptstyle
e_3$}}}}
\put(1011,3420){\makebox(0,0)[lb]{\smash{{\SetFigFont{12}{14.4}{\rmdefault}{\mddefault}{\updefault}$\scriptstyle
e_4$}}}}
\put(2026,5130){\makebox(0,0)[lb]{\smash{{\SetFigFont{12}{14.4}{\rmdefault}{\mddefault}{\updefault}$\scriptstyle
e_6$}}}}
\put(2076,2930){\makebox(0,0)[lb]{\smash{{\SetFigFont{12}{14.4}{\rmdefault}{\mddefault}{\updefault}$\scriptstyle
e_1$}}}}
\put(3046,3420){\makebox(0,0)[lb]{\smash{{\SetFigFont{12}{14.4}{\rmdefault}{\mddefault}{\updefault}$\scriptstyle
e_1$}}}}
\put(2081,1780){\makebox(0,0)[lb]{\smash{{\SetFigFont{12}{14.4}{\rmdefault}{\mddefault}{\updefault}$\scriptstyle
e_2$}}}}
\put(4051,1770){\makebox(0,0)[lb]{\smash{{\SetFigFont{12}{14.4}{\rmdefault}{\mddefault}{\updefault}$\scriptstyle
e_5$}}}}
\end{picture}
}
\begin{center}
Figure 1
\end{center}
\end{figure}



Every octonion $q\in \OO$ can be written uniquely in the form
$$q=\sum_{i=0}^7x_ie_i$$
where $x_i\in \RR$. The summand $x_0e_0=x_0$ is called the real
part of $q$ and is denoted by $Re(q)$.

One defines the octonionic conjugate of $q$ by
$$\bar q:=x_0-\sum_{i=1}^7x_ie_i.$$
It is well known that the conjugation is an anti-involution of
$\OO$:
$$\overline{\overline{q}}=q,\, \overline{a+b}=\bar a+\bar b,\, \overline{ab}=\bar b\bar a.$$
Let us define a norm on $\OO$ by $|q|:=\sqrt{q\bar q}$. Then
$|\cdot|$ is a multiplicative norm on $\OO$: $|ab|=|a||b|$. The
square of the norm $|\cdot|^2$ is a positive definite quadratic
form. Its polarization is a  positive definite scalar product
$<\cdot,\cdot>$ on $\OO$ which is given explicitly by
$$<x,y>=Re(x\bar y).$$

Furthermore $\OO$ is a division algebra: any $q\ne 0$ has a unique
inverse $q^{-1}$ such that $qq^{-1}=q^{-1}q=1$. In fact
$$q^{-1}=|q|^{-2}\bar q.$$

We denote by $\HH$ the usual quaternions. It is associative division
algebra. We will fix once and for all an imbedding of algebras
$\HH\subset \OO$. Let us denote by $i,j\in \HH$ the usual
quaternionic units, and $k=ij$. Then $i,j,k$ are pairwise orthogonal
with respect to the scalar product $<\cdot,\cdot>$. Let us fix once
and for all an octonionic unit $l\in \OO, l^2=-1$ which is
orthogonal to $i,j,k$. Then $l$ anti-commutes with $i,j,k$. Every
element $q\in \OO$ can be written uniquely in the form
$$q=x+yl \mbox{ with } x,y\in \HH.$$
Then we can multiply two such octonions using the formula
\begin{eqnarray}\label{E:product-oct}
(x+yl)(w+zl)=(xw-\bar z y)+(zx+y\bar w)l
\end{eqnarray}
where $x,y,w,z\in\HH$.

We have the following weak forms of the associativity in
octonions.
\begin{lemma}\label{L:oalg1}
Let $a,b,c\in \OO$. Then

(i) $Re((ab)c)=Re(a(bc))$ (this real number will be denoted by
$Re(abc)$).

(ii) $a(bc)+\bar b(\bar a c)=(ab+\bar b\bar a)c$.

(iii) $(ca)b+(c\bar b)\bar a=c(ab+\bar b \bar a)$.

(iv) Any subalgebra of $\OO$ generated by any two elements and
their conjugates is associative. It is always isomorphic either to
$\RR$, $\CC$, or $\HH$.

(v) $Re((\bar ab)( ca))=|a|^2Re(bc).$
\end{lemma}
{\bf Proof.} For (i) see Corollary 15.12(i) in \cite{adams-book}.
For (iv) see e.g. Chapter 15 in \cite{adams-book}, particularly
Lemma 15.6 which is essentially equivalent to statement (iv) of
our lemma.

Observe that (iii) is equivalent to (ii) by taking the
conjugation. Thus let us prove (ii). Let us defines the 3-linear
map $[\cdot,\cdot,\cdot]\colon\OO\times\OO\times \OO\to \OO$
called associator which is defined by $[x,y,z]:=(xy)z-x(yz)$. Then
the statement (ii) is equivalent to
\begin{eqnarray}\label{E:assoc}
[a,b,c]=-[\bar b,\bar a,c].
\end{eqnarray}
By Theorem 15.11(ii) of \cite{adams-book} the associator $[x,y,z]$
changes sign when one conjugates any variables. By Theorem
15.11(iii) of \cite{adams-book} the associator $[x,y,z]$ is an
alternating function of three variables. These two properties
imply the identity (\ref{E:assoc}).

Let us prove part (v). It is clear that both sides of the equality
we have to prove, are linear in $b,c$. By the part (iv) they are
equal to each other (even without taking the real part) in the
following two cases: (i) at least one of $b$ and $c$ is real; (ii)
$b$ and $c$ are proportional to each other (with a real
coefficient). Hence we may assume that $Re (b)=Re (c)=0$ and
$<b,c>=0$. By applying an appropriate automorphism of $\OO$ we may
assume that $b=i,c=j$. Then it remains to show that for any
$a\in\OO$
$$Re((\bar ai)( ja))=0.$$
This can be check by a direct computation using a representation
$a=u+vl$ with $u,v\in\HH$ and the formula (\ref{E:product-oct}).
Lemma is proved. \qed

\subsection{Octonionic hermitian matrices.}\label{S:la}\setcounter{theorem}{0}
Let us denote
$$\OO^m:=\{(q_1,\dots,q_m)|\, q_i\in \OO\}.$$
For $\xi=(q_1,\dots,q_m)\in \OO^m$, $x\in \OO$ we will denote by
$\xi\cdot x$ the $m$-tuple $(q_1x,\dots,q_mx)\in\OO^m$. Notice that
usually we will write elements of $\OO^m$ as $m$-columns rather than
rows.

Let us denote by $\ch_n(\RR)$ the space of real symmetric
$(n\times n)$-matrices. The space $\ch_n(\RR)$ is naturally
identified with the space of real valued quadratic forms on
$\RR^n$.

Let us denote by $\ch_m(\OO)$ the space of octonionic hermitian
$(m\times m)$-matrices. By definition, an $(m\times m)$-matrix
$A=(a_{ij})$ with octonionic entries is called hermitian if
$a_{ij}=\overline{a_{ji}}$ for any $i,j$. For a matric $A=(a_{ij})$
denote also $A^*:=(\overline{a_{ji}})$. We will be mostly interested
in octonionic hermitian matrices of size 2. In this case we have the
following explicit description
\begin{eqnarray}
\ch_2(\OO)=\left\{ \left[\begin{array}{cc}
                          a     &q\\
                          \bar q&b
                         \end{array}\right]\big| a,b\in \RR, q\in \OO\right\}.
\end{eqnarray}

We have the natural $\RR$-linear map
\begin{eqnarray}\label{D:j}
j\colon\ch_2(\OO)\to \ch_{16}(\RR)
\end{eqnarray} which is defined
as follows: for any $A\in \ch_2(\OO)$ the value of the quadratic
form $j(A)$ on any octonionic 2-column $\xi\in \OO^2\simeq
\RR^{16}$ is equal $j(A)(\xi)=Re(\xi^* A\xi)$ (note that the
bracketing inside the formula is not important due to Lemma
\ref{L:oalg1}(i)). It is easy to see that the map $j$ is
injective. Via this map $j$ we will identify $\ch_2(\RR)$ with a
subspace of $\ch_{16}(\RR)$.

\def\hr{\ch_{16}(\RR)}
\def\ho{\ch_2(\OO)}
Let us construct now a linear map
$$\theta\colon \ch_{16}(\RR)\to \ch_2(\OO)$$
such that $\theta\circ j=Id$ and which will be useful later. For any
$B\in \ch_{16}(\RR)$ let us denote by $b$ the corresponding
quadratic form on $\RR^{16}\simeq \OO^2$. Define
$$\theta(B):=\frac{1}{16}\left(\ddbq\right)_{i,j=1}^2$$
to be the octonionic Hessian of $b$ all of whose entries are
replaced by octonionic conjugates. Note that the matrix
in the right hand side of the last formula is independent of a point
in $\OO^2$.
\begin{lemma}\label{L:la1}
For any $B\in \hr$ and any $\xi\in \OO^2$ of the form either
$\xi=\left[\begin{array}{c}
                                          a\\
                                          1
                                          \end{array}\right]$ or
$\xi=\left[\begin{array}{c}
                                          1\\
                                          a
                                          \end{array}\right]$, $a\in \OO$, one has
$$\theta(B)(\xi):=Re\left(\xi^*\theta(B)\xi\right)= \int_{x\in \OO,|x|=1}b(\xi\cdot x)dx$$
where $dx$ is the rotation invariant probability measure on
$S^7=\{x\in \OO,|x|=1\}$.
\end{lemma}
{\bf Proof.} Let us denote elements of $\RR^{16}$ by 16-tuples
$(x_0,x_1,\dots,x_7;y_0,y_1,\dots,y_7)$ where $x_i$'s correspond to
the first octonionic coordinate $q_1=\sum_{i=0}^7x_ie_i$, and
similarly $y_i$'s correspond to the second quaternionic coordinate
$q_2$. By linearity and symmetry considerations it is enough to
prove the lemma in the following 2 cases for fixed $p,q$:

(1) $b((x,y))=x_px_q$;

(2) $b((x,y))=x_py_q$ for any $(x,y)\in \OO^2$.

Let us start with case (1). For any $\xi=\left[\begin{array}{c}
                                          a\\
                                          1
                                          \end{array}\right]\in \OO^2$ one has
\begin{eqnarray*}
\int_{x\in \OO,|x|=1}b(\xi\cdot
x)dx=\int_{S^7}(ax)_p(ax)_qdx=|a|^2\int_{S^7}x_px_qdx.
\end{eqnarray*}
The last integral vanishes for $p\ne q$ and equals to
$\frac{|a|^2}{8}$ for $p=q$. On the other hand if $p\ne q$
$\left(\ddbq\right)=0$, thus lemma is proved in this case. If $p=q$
then
$$\left(\ddbq\right)=\left[\begin{array}{cc}
              2&0\\
              0&0
             \end{array}\right].$$

Hence $$\theta(B)(\xi)=\frac{1}{16} Re\left([\bar a,
1]\left[\begin{array}{cc}
              2&0\\
              0&0
             \end{array}
\right]\left[\begin{array}{c}
                                         a\\
                                         1
                                        \end{array}\right]\right)=
\frac{|a|^2}{8}=\int_{x\in \OO,|x|=1}b(\xi\cdot x)dx.$$

Let us consider case (2). In this case
$$\left(\ddbq\right)=\left[\begin{array}{cc}
                          0&e_p\bar e_q\\
                          e_q\bar e_p&0
                          \end{array}\right].$$
Let first $\xi=\left[\begin{array}{c}
                                          a\\
                                          1
                                          \end{array}\right]$. Then
$$\theta(B)(\xi)=\frac{1}{16}Re\left([\bar a,1]\left[\begin{array}{cc}
                          0&e_p\bar e_q\\
                          e_q\bar e_p&0
                          \end{array}\right]\left[\begin{array}{c}
                                             a\\
                                             1
                                             \end{array}\right]\right)=
\frac{1}{8}Re(e_q\bar e_p a).$$ On the other hand
\begin{eqnarray*}
\int_{x\in \OO,|x|=1}b(\xi\cdot x)d(x)=
\int_{S^7}(ax)_px_qd(x)=\int_{S^7} Re(\bar e_p ax)\cdot x_q dx
=\\\int_{S^7} Re\left(\bar e_p a(\sum_{s=0}^7x_se_s)\right)x_q
dx=\int_{S^7}Re(\bar e_p a e_q)x_q^2=\frac{1}{8}Re(e_q\bar e_p
a)=\theta(B)(\xi).
\end{eqnarray*} Similarly one considers the case
$\xi=\left[\begin{array}{c}
                                          1\\
                                          a
                                          \end{array}\right]$.
Lemma is proved. \qed

\begin{corollary}\label{C:la2}
$$\theta\circ j=Id.$$
\end{corollary}
{\bf Proof.} It is easy to see that any matrix
$A=\left[\begin{array}{cc}
          \alp&q\\
          \bar q&\beta
          \end{array}\right] \in \ho$ is uniquely determined by the
products $Re(\xi^*A\xi)$ where $\xi=\left[\begin{array}{c}
                                          a\\
                                          1
                                          \end{array}\right]\in \OO^2$. By Lemma \ref{L:la1} one has
\begin{eqnarray*}\left((\theta\circ
j)(A)\right)(\xi)=\int_{S^7}(jA)(\xi\cdot x)dx= \int_{S^7}
Re\left([\overline{ax},\bar x]\left[\begin{array}{cc}
          \alp&q\\
          \bar q&\beta
          \end{array}\right]\left[\begin{array}{c}
                                   ax\\
                                   x
                                   \end{array}\right]\right)dx=\\
\int_{S^7}(\alp|a|^2+\beta+2Re((\overline{ax})qx))dx\overset{\mbox{Lemma
}\ref{L:oalg1}(v)}{=} \alp|a|^2+\beta+2Re(q\bar
a)=\\Re\left([\overline{a}, 1]\left[\begin{array}{cc}
          \alp&q\\
          \bar q&\beta
          \end{array}\right]\left[\begin{array}{c}
                                   a\\
                                   1
                                   \end{array}\right]\right)=Re(\xi^*A\xi).
\end{eqnarray*}
Corollary is proved. \qed

\begin{definition}\label{D:positive-def}
Let $A\in\ch_2(\OO)$. $A$ is called {\itshape positive definite}
(resp. {\itshape non-negative definite}) if for any 2-column
$\xi\in\OO^2\backslash\{0\}$
$$Re(\xi^* A\xi)>0 \,\, (\mbox{resp. } Re(\xi^* A\xi)\geq 0).$$
\end{definition}
For a positive definite (resp. non-negative definite) matrix $A$ one
writes as usual $A>0$ (resp. $A\geq 0$).

On the class of octonionic hermitian $(2\times 2)$-matrices there is
a nice notion of determinant which is defined by
\begin{eqnarray}\label{D:det}
\det\left(\left[\begin{array}{cc}
                a&q\\
                \bar q&b
                \end{array}\right]\right)=ab-|q|^2.
\end{eqnarray}
\begin{remark}
It turns our that a nice notion of determinant does exist also on
octonionic hermitian matrices of size 3, see e.g. Section 3.4 of
\cite{baez}. Note also that a nice determinant does exist for
{\itshape quaternionic} hermitian matrices of any size: see the
survey \cite{aslaksen}, the article \cite{gelfand-retakh-wilson},
and for applications to quaternionic plurisubharmonic functions see
\cite{alesker-busm-03}, \cite{alesker-jga-03},
\cite{alesker-adv-05}, \cite{alesker-alg-anal}.
\end{remark}

The following result is a version of the Sylvester criterion for
octonionic matrices of size two.
\begin{proposition}\label{P:sylvester}
Let $A=\left[\begin{array}{cc}
              a&q\\
              \bar q&b
              \end{array}\right]\in\ch_2(\OO)$. Then $A>0$ if and
only if $a>0$ and $\det A>0$.
\end{proposition}
{\bf Proof.} Let $\xi=\left[\begin{array}{c}
                             x\\
                             y
                             \end{array}\right]\in \ch_2(\OO)$.
Then
\begin{eqnarray}\label{E:ca-sh}
Re(\xi^*A\xi)=a|x|^2+b|y|^2+ 2Re(qy\bar x).
\end{eqnarray}
Assume first that $A>0$. Taking $x=1,y=0$ we get $a>0$. Next
substituting into (\ref{E:ca-sh}) $y=t|q|q^{-1} x$ where $t\in \RR$
we get for any $x\in \OO\backslash\{0\}$ and any $t\in \RR$
$$|x|^2(a+bt^2+2t|q|) >0.$$ hence $0<ab-|q|^2=\det A$.

Conversely assume that $a>0,\det A>0$. Since $|Re(uv)|\leq |u|\cdot
|u|$, (\ref{E:ca-sh}) implies that
$$Re(\xi^*A\xi)\geq a|x|^2+b|y|^2-2|q|\cdot |x|\cdot|y|.$$
Our assumptions imply that the last expression is positive
provided $(x,y)\ne 0$. \qed

Now will will introduce the notion of {\itshape mixed determinant}
of two octonionic hermitian matrices in analogy to the classical
real case (see e.g. \cite{schneider-book}). First observe that the
determinant $\det \ch_2(\OO)\to \RR$ is a homogeneous polynomial of
degree 2 on the real vector space $\ch_2(\OO)$. Hence it admits a
unique polarization: a bilinear symmetric map
$$D\colon \ch_2(\OO)\times \ch_2(\OO)\to \RR$$
such that $D(A,A)=\det A$ for any $A\in \ch_2(\OO)$. This map $D$ is
called the mixed determinant. By the abuse of notation it will be
denoted again by $\det$. Explicitly if $A=(a_{ij})_{i,j=1}^2$,
$B=(b_{ij})_{i,j=1}^2$ are octonionic hermitian then
\begin{eqnarray}\label{E:mixed-det}
\det(A,B)=\frac{1}{2}\left(a_{11}b_{22}+a_{22}b_{11}-2Re(a_{12}b_{21})\right).
\end{eqnarray}
\begin{lemma}\label{L:mix-det-posit}
If $A,B\in \ch_2(\OO)$ are positive (resp. non-negative) definite
then
$$\det(A,B)>0 (\mbox{ resp. } \det(A,B)\geq 0).$$
\end{lemma}
{\bf Proof.} Let us assume that $A,B$ are positive definite. Then by
Proposition \ref{P:sylvester} we have
\begin{eqnarray*}
a_{11}>0,a_{22}>0, |a_{12}|<\sqrt{a_{11}a_{22}},\\
b_{11}>0,b_{22}>0, |b_{12}|<\sqrt{b_{11}b_{22}}.
\end{eqnarray*}
These inequalities imply
\begin{eqnarray}\label{E:iii}
Re(a_{12}b_{21})\leq |a_{12}|\cdot
|b_{21}|<\sqrt{a_{11}a_{22}b_{11}b_{22}}.
\end{eqnarray} Substituting
(\ref{E:iii}) into (\ref{E:mixed-det}) we get
\begin{eqnarray*}
\det(A,B)>\frac{1}{2}\left(a_{11}b_{22}+a_{22}b_{11}-
2\sqrt{a_{11}a_{22}b_{11}b_{22}}\right)\geq 0
\end{eqnarray*}
where the last estimate is the arithmetic-geometric mean inequality.
Thus $\det(A,B)>0$ for positive definite matrices $A,B$. For
non-negative definite matrices the result follows by going to the
limit. \qed

\begin{remark}
Also one can prove the following version of the Aleksandrov
inequality for mixed determinants \cite{aleksandrov-mixed}: The
mixed determinant $\det(\cdot,\cdot)$ is a non-degenerate quadratic
form on the real vector space $\ch_2(\OO)$; its signature type has
one plus and the rest are minuses. Consequently, if $A>0$ then for
any $B\in \ch_2(\OO)$ one has
$$\det(A,B)^2\geq \det A\cdot \det B$$
and the equality is achieved if and only if $A$ is proportional to
$B$ with a real coefficient.

We do not present a detailed proof of this result since we are not
going to use it. Notice only that this result can be deduced
formally from the corresponding result for quaternionic matrices
proved in \cite{alesker-busm-03}. Indeed all the entries of $A$ and
$B$ together  contain at most two non-real octonions, hence the
field generated them is associative by Lemma \ref{L:oalg1}(iv).
\end{remark}

\subsection{Octonionic projective line
$\OO\PP^1$.}\label{Ss:projline} In this section we remind the
definition and basic properties of the octonionic projective line
$\OO\PP^1$.

Let us define an equivalence relation $\sim$ on the unit sphere
$S^{15}\subset \OO^2$ by saying $\xi\sim \eta$ if and only if
$$\xi\xi^*=\eta\eta^*$$
where we write $\xi$ and $\eta$ as columns, thus $\xi\xi^*$ and
$\eta\eta^*$ are octonionic $(2\times 2)$-matrices. It is easy to
see that $\left[\begin{array}{c}
                 x\\
                 y
                 \end{array}\right]\sim \left[\begin{array}{c}
                                                xy^{-1}\\
                                                1
                                                \end{array}\right]$
                                                if $y\ne 0$, and
$\left[\begin{array}{c}
                 x\\
                 y
                 \end{array}\right]\sim \left[\begin{array}{c}
                                                1\\
                                                yx^{-1}
                                                \end{array}\right]$
                                                if $x\ne 0$.
The quotient of $S^{15}$ by this equivalence relation is called the
octonionic projective line $\OO\PP^1$.
\begin{remark}
If in the above construction one replaces $\OO$ by $\RR,\CC,$ or
$\HH$, one get the usual projective lines
$\RR\PP^1,\CC\PP^1,\HH\PP^1$ respectively.
\end{remark}
$\OO\PP^1$ has a natural smooth structure, and it is diffeomorphic
to the standard sphere $S^8$. The fibers of the quotient map
$S^{15}\to \OO\PP^1$ are the spheres $S^7$. This map is called the
octonionic Hopf fibration.

\subsection{The group $SL_2(\OO)$.}\label{Ss:sl2}
We discuss in this section the definition and basic properties of
the group $SL_2(\OO)$. We refer to \cite{sudbery} for the proofs and
further details.

An octonionic $(2\times 2)$-matrix is called {\itshape traceless} if
the sum of its diagonal elements is equal to zero. Every $(2\times
2)$-matrix $A$ with octonionic entries defines an $\RR$-linear
operator on $\OO^2$ by $\xi\mapsto A\cdot\xi$. However the space of
such operators is not closed under the commutator due to the lack of
associativity. One denotes by $sl_2(\OO)$ the Lie subalgebra of
$gl_{16}(\RR)$ generated by $\RR$-linear operators on
$\OO^2\simeq \RR^{16}$ determined by all traceless octonionic
matrices. This Lie algebra $sl_2(\OO)$ turns our to be semi-simple
\cite{sudbery} (see Theorem \ref{T:sl2-descr} below for details).
But any semi-simple Lie subalgebra of an algebraic group is a Lie
algebra of a closed algebraic subgroup (see e.g.
\cite{onishik-vinberg}, Ch. 3, \S 3.3). In our case this subgroup is
denoted by $SL_2(\OO)\subset GL(16,\RR)$.

\begin{theorem}[\cite{sudbery}]\label{T:sl2-descr}
(i) The Lie algebra $sl_2(\OO)$ is isomorphic to the Lie algebra
$so(9,1)$.

(ii) The Lie group $SL_2(\OO)$ is isomorphic to the group
$Spin(9,1)$ (which is the universal covering of the identity
component of the pseudo-orthogonal group $O(9,1)$).
\end{theorem}
\begin{remark}\label{R:spin9}
A maximal compact subgroup of $SL_2(\OO)\simeq Spin(9,1)$ is
isomorphic to the group $Spin(9)$ which is the universal covering of
the special orthogonal group $SO(9)$.
\end{remark}

Both $sl_2(\OO)$ and $SL_2(\OO)$ come with their {\itshape
fundamental representations} on $\OO^2$. Moreover the Lie algebra
$sl_2(\OO)$ acts on the space $\ch_2(\OO)$. This action is uniquely
characterized by the following property (see \cite{sudbery} for the
details): if $A$ is a traceless matrix it acts by $A\colon X\mapsto
-A^*X-XA$. Since the group $SL_2(\OO)$ is connected and simply
connected, this representation of $sl_2(\OO)$ integrates to a
representation of the group $SL_2(\OO)$ on $\ch_2(\OO)$.

\begin{proposition}\label{C:cone-preserved}
The group $SL_2(\OO)$ preserves the cone of positive definite
octonionic hermitian matrices.
\end{proposition}
\def\clk{\overline{\ck}}
{\bf Proof.} Let us denote by $\ck$ the open cone of positive
definite matrices in $\ch_2(\OO)$. Let us denote by $\overline{\ck}$
the closure of $\ck$, namely the closed cone of non-negative
definite matrices. The boundary $\pt\bar\ck$ is a hypersurface in
$\ch_2(\OO)$ which is smooth at every point except of $0$.

In order to prove the proposition it is enough to prove the
infinitesimal version of it as follows. Let us consider any element
$D\in sl_2(\OO)$ from the Lie algebra. It induces a vector field on
$\ch_2(\OO)$ via its action: $X\mapsto D(X)$. In order to show that
the one-parametric subgroup in $SL_2(\OO)$ generated by $D$
preserves the cone $\ck$ it is enough to check that at any point
$X\in \pt\bar\ck$ the vector $D(X)$ is not directed outside of the
domain $\clk$ (i.e. looks inside or tangent to the boundary
$\pt\clk$) when $X$ is a smooth point of the boundary, and vanishes
when $X$ is a singular point of the boundary. Clearly $D(0)=0$, and
$0$ is the only singular point of the boundary. Hence we may assume
that $X$ is a smooth point of the boundary. We are going to show
that the vector $D(X)$ is in fact tangent to $\pt\clk$. Since
$\pt\clk=\{A\geq 0|\, \det A=0\},$ this follows from the fact that
the group $SL_2(\OO)$ preserves the determinant of octonionic
hermitian matrices (see e.g. \cite{baez}, p. 177). Proposition is
proved. \qed

The following two lemmas are essentially contained in
\cite{manogue-schray}, but we would like to present a proof for the
sake of completeness.
\begin{lemma}[\cite{manogue-schray}]\label{m-sch1}
The linear map $\OO^2\otimes_\RR \OO^2\to \ch_2(\OO)$ given by
$$\xi\otimes \eta\mapsto \xi\cdot \eta^*+\eta\cdot\xi^*$$
is $SL_2(\OO)$-equivariant.
\end{lemma}
{\bf Proof.} It is enough to prove the equivariance with respect to
the action of the Lie algebra $sl_2(\OO)$. In fact it is enough to
show that for any octonionic traceless $(2\times 2)$-matrix $M$ and
any $\xi\in \OO^2$ one has
\begin{eqnarray}\label{ff1}
(M\xi)\cdot \xi^*+\xi\cdot (M\xi)^*=M\cdot
(\xi\xi^*)+(\xi\xi^*)\cdot M^*.
\end{eqnarray}
This equality is proved by a straightforward computation using Lemma
\ref{L:oalg1}. \qed

\begin{lemma}[\cite{manogue-schray}]\label{m-sch2}
(i) The group $SL_2(\OO)$ acts naturally on $\OO\PP^1$, namely for
any $\phi\in SL_2(\OO)$ and any $L\in \OO\PP^1$ the subspace
$\phi(L)$ is an octonionic projective line.

(ii) For any $L\in \OO\PP^1$ and any $\phi\in SL_2(\OO)$ the
restriction
$$\phi|_L\colon L\to \phi(L)$$
is a conformal linear map.
\end{lemma}
{\bf Proof.} First let us show that if $\xi,\eta\in \OO^2$ have norm
1 and $\xi\sim \eta$ then $|\phi(\xi)|=|\phi(\eta)|$. Observe that
for any $v\in \OO^2$ one has
$$|v|^2=v^*v=Tr(vv^*)$$
where $Tr$ denotes the sum of diagonal elements of a matrix. Then we
have
\begin{eqnarray*}
|\phi(\xi)|^2=Tr\left(\phi(\xi)\cdot
\phi(\xi)^*\right)\overset{\mbox{by Lemma
}\ref{m-sch1}}{=}Tr(\phi(\xi\cdot
\xi^*))=\\Tr(\phi(\eta\cdot\eta^*))\overset{\mbox{by Lemma
}\ref{m-sch1}}{=}Tr\left(\phi(\eta)\cdot
\phi(\eta)^*\right)=|\phi(\eta)|^2.
\end{eqnarray*}
Thus $ |\phi(\xi)|=|\phi(\eta)|$.Then in order to prove both parts
of the lemma it remains to show that for any norm 1 vectors
$\xi,\eta\in \OO^2$ such that $\xi\sim \eta$ (i.e.
$\xi\xi^*=\eta\eta^*$) one has $\phi(\xi)\sim \phi(\eta)$, i.e.
$\phi(\xi)\cdot \phi(\xi)^*=\phi(\eta)\cdot \phi(\eta^*)$. But
applying Lemma \ref{m-sch1} twice we get
$$\phi(\xi)\cdot \phi(\xi)^*=
\phi(\xi\xi^*)=\phi(\eta\eta^*)=\phi(\eta)\cdot\phi(\eta)^*.$$ \qed

We immediately deduce the following corollary.
\begin{corollary}
The octonionic Hopf map $S^{15}\to \OO\PP^1$ is
$Spin(9)$-equivariant.
\end{corollary}

The following lemma is well known \cite{borel2}.
\begin{lemma}
The group $Spin(9)$ acts transitively on the unit sphere $S^{15}$,
and hence on $\OO\PP^1$.
\end{lemma}

\subsection{Further properties of the octonionic
Hessian.}\label{Ss:octHessian} \setcounter{theorem}{0}
\begin{proposition}\label{P:ohs1}
Let $f\colon \OO^2\to \RR$ be a $C^2$-smooth function. Let $A\in
SL_2(\OO)$. Then
$$\left(\frac{\pt^2}{\pt\bar q_i \pt 
q_j}(f(A^{-1}q))\right)=A\left(\dfq (A^{-1}q)\right)$$ where $A$ in the right
hand side denotes the induced action of $A$ on $\ho$.
\end{proposition}
{\bf Proof.} By translation it is enough to check the above
equality at $q=0$. Moreover we may and will assume that $f$ is a
quadratic form. Thus the proposition becomes equivalent to the the
following lemma.

\begin{lemma}\label{L:theta-equivar}
The map $\theta\colon \hr\to \ho$ is $SL_2(\OO)$-equivariant.
\end{lemma}
{\bf Proof.} It is enough to prove this proposition infinitesimally,
i.e. for the action of the Lie algebra $sl_2(\OO)$. Moreover it is
sufficient to check it for a set of generators, say for traceless
$(2\times 2)$-octonionic matrices. Let us fix such a traceless
matrix $A$. We have to show that for any $b\in \ch_{16}(\RR)$ one
has
\begin{eqnarray}
\label{bb1}\theta(A(b))=A(\theta(b)).
\end{eqnarray}
It is enough to show that for any $\xi=\left[\begin{array}{c}
                                                   p\\
                                                   1\end{array}\right]\in \OO^2$
one has
\begin{eqnarray}\label{bb1.5}
Re(\xi^*\theta(A(b))\xi)=Re(\xi^*A(\theta(b))\xi).
\end{eqnarray}
The right hand side of (\ref{bb1.5}) can be
rewritten as
\begin{eqnarray}\label{bb1.8}
-Re\left(\xi^*(A^*\theta(b))\xi+\xi^*(\theta(b)A)\xi\right)=
-2Re\left(\xi^*(A^*\theta(b))\xi\right).
\end{eqnarray}

By Lemma \ref{L:la1} the left hand side of (\ref{bb1.5}) can be
rewritten as follows:
\begin{eqnarray}\label{bb2}
-2\int_{x\in S^7}b((A\xi)x,\xi
x)dx=-\frac{d}{d\tau}\big|_0\int_{x\in S^7}b((A\xi+\tau\xi)x)dx=
\\-\frac{d}{d\tau}\big|_0\theta(b)(A\xi+\tau\xi)=
-\frac{d}{d\tau}\big|_0Re\left((A\xi+\tau\xi)^*\theta(b)(A\xi+\tau\xi)\right)=\\\label{bb3}
-2Re\left((\xi^*A^*)\theta(b)\xi\right).
\end{eqnarray}
Substituting (\ref{bb1.8}) and (\ref{bb3}) into (\ref{bb1.5}) we see
that (\ref{bb1.5}) becomes equivalent to
\begin{eqnarray}\label{bb4}
Re\left(\xi^*(A^*\theta(b))\xi\right)=Re\left((\xi^*A^*)\theta(b)\xi\right).
\end{eqnarray}
Both sides of the above quality are linear with respect to
$\theta(b)\in\ch_2(\OO)$. Obviously the equality is satisfied when
$\theta(b)$ has real entries. Thus we may assume that
$\theta(b)=\left[\begin{array}{cc}
                      0&q\\
                      \bar q&0
                      \end{array}\right]$. Let us denote also
$A^*=\left[\begin{array}{cc}
            a&b\\
            c&d
            \end{array}\right]$. In this notation (\ref{bb4})
becomes
\begin{eqnarray}\label{bb5}
Re\left([\bar p,1]\left(\left[\begin{array}{cc}
                               a&b\\
                               c&d
                               \end{array}\right]\left[\begin{array}{cc}
                                                   0&q\\
                                                   \bar q&0
                                                   \end{array}\right]\right)\left[\begin{array}{c}
                                                                 p\\
                                                                 1
                                                                 \end{array}\right]\right)=
Re\left(\left([\bar p,1]\left[\begin{array}{cc}
                               a&b\\
                               c&d
                               \end{array}\right]\right)\left[\begin{array}{cc}
                                                   0&q\\
                                                   \bar q&0
                                                   \end{array}\right]\left[\begin{array}{c}
                                                                 p\\
                                                                 1
                                                                 \end{array}\right]\right)
\end{eqnarray}

By a direct computation the right hand side of (\ref{bb5}) is equal
to
\begin{eqnarray}\label{bb6}
Re((\bar pb)(\bar qp))+Re\left(\bar p aq+cq+d\bar qp\right).
\end{eqnarray}
The left hand side of (\ref{bb5}) is equal to
\begin{eqnarray}\label{bb7}
|p|^2Re(b\bar q)+Re\left( \bar p aq+cq+d\bar qp\right).
\end{eqnarray}


Comparing (\ref{bb6}) and (\ref{bb7}), it remains to prove that for
any $b,p,q\in\OO$ one has
\begin{eqnarray}\label{bb8}
Re((\bar pb)( qp))=|p|^2Re(b q).
\end{eqnarray}
But this is exactly Lemma \ref{L:oalg1}. Hence Lemma
\ref{L:theta-equivar} is proved. \qed


\begin{remark}
Similarly one can show that the imbedding $j\colon
\ch_2(\OO)\inj\ch_{16}(\RR)$ (see (\ref{D:j}) in Section \ref{S:la})
is also $SL_2(\OO)$-equivariant.
\end{remark}


\section{Octonionic Radon transform.}\label{S:radon}
\setcounter{theorem}{0}\setcounter{subsection}{1}
\begin{definition}\label{D:radon1}
An {\itshape affine octonionic line} in $\OO^2$ is any translate
of an octonionic line from $\OO\PP^1$.
\end{definition}
The manifold of all affine octonionic lines in $\OO^2$ will be
denoted by $\ca\OO\PP^1$. It is a homogeneous space for the group
$\OO^2\rtimes Spin(9)$ (we denote in this way the semi-direct
product of the group $\OO^2$ of translations of $\OO^2$, and of the
group $Spin(9)$ of certain linear transformations).

Let us define the Radon transform operator
\begin{eqnarray}\label{radon1.5}
R\colon C^\infty_0(\OO^2)\to C^\infty_0(\ca\OO\PP^1)
\end{eqnarray}
by $(Rf)(E)=\int_{q\in E}f(q)dq$ where $dq$ is the Lebesgue
measure induced by the standard Euclidean metric on $\OO^2$.
\begin{proposition}\label{P:radon2}
The octonionic Radon transform (\ref{radon1.5}) is injective.
\end{proposition}
{\bf Proof.}  We will just present the inversion formula completely
analogous to the complex Radon transform (see \cite{gelfand}; see
also appendix in \cite{alesker-adv-05} for the quaternionic case).
For any point $q\in \OO^2$ let ${\cal P}_q$ denote the manifold of
affine octonionic lines passing through $q$. Clearly $\cp_q\simeq
\OO\PP^1$. For $E\in \ca\OO\PP^1$ let us denote by $E^{\perp}$ the
octonionic line orthogonal to $E$ and passing through the origin 0.

Let us define the  operator
$$\cd: C^{\infty}(\ca\OO\PP^1) \to C^{\infty}(\OO ^2) $$
as follows. Let $g\in C^{\infty}(\ca\OO\PP^1) $. Set  $$\cd g(q):=
\int_{E\in {\cp}_q} ( \Delta _{E^\perp}) ^{4} g(E+w)dE,$$ where
$\Delta _{E^\perp}$ denotes the (8-dimensional) Laplacian with
respect to $w\in E^\perp$, and the integration is with respect the
Haar measure on ${\cal P}_q$ invariant under the action of
$Spin(9)$.

\begin{claim}
For any smooth rapidly decreasing function $f$ of $\OO^2$
$$\cd(Rf)= c\cdot f,$$
where $c$ is a non-zero constant.\end{claim}

It is sufficient to check this claim pointwise, say at $0$. The
operators $R$ and $\cd$ commute with translations and the action of
the group $Spin(9)$. Then $\cd(Rf)(0)$ defines a distribution which
is invariant with respect to the action of $Spin(9)$. Moreover it is
easy to check that this distribution is homogeneous of degree $-16$
(exactly as the delta-function at $0$). Since the group $Spin(9)$
acts transitively on the unit sphere $S^{15}$, there is at most one
dimensional space of $Spin(9)$- invariant distributions homogeneous
of degree $-16$. Hence they must be proportional to the
delta-function at $0$. Thus $\cd(Rf)= c\cdot f$ for some constant
$c$. To see that $c\ne 0$ it is sufficient to check it by an
explicit computation for the function $f(q)=exp(-|q|^2/2)$. \qed


For an affine line $L\in \ca\OO\PP^1$ let us denote by $\delta_L$
the generalized function given by
$$\delta_L(\mu)=\int_L\mu$$
for any infinitely smooth compactly supported measure $\mu$.

\begin{corollary}\label{C:radon3}
The $\CC$-linear span of $\delta$-functions of affine octonionic
lines is dense in the weak topology in the space
$C^{-\infty}(\OO^2)$ of generalized functions.
\end{corollary}
{\bf Proof.} This follows immediately from Proposition
\ref{P:radon2} and the Hahn-Banach theorem. \qed

\section{Octonionic plurisubharmonic functions and their properties.}\label{S:psh}

\subsection{Octonionic plurisubharmonic
functions.}\label{Ss:psh-func} Let us remind few standard
definitions. Recall that a real valued function $f$ is called upper
semi-continuous if $$f(x_0)\geq \lim\sup_{x\to x_0}f(x)$$ for any
point $x_0$.

\begin{definition}\label{D:sh}
Let $U\subset \RR^n$ be an open subset. A function $f\colon U\to
\RR\cup\{-\infty\}$ is called {\itshape subharmonic} if

(i) $f$ is upper semi-continuous;

(ii) for any point $x_0\in U$ and any sphere $S\subset U$ centered
at $x_0$ one has
\begin{eqnarray}\label{E:sh}
f(x_0)\leq\int_Sf(x)dx
\end{eqnarray}
where $dx$ is the probability rotation invariant measure on the
sphere $S$.
\end{definition}
According to this definition the function which is identically equal
to $-\infty$ is subharmonic.
\begin{remark}
(1) The integral in (\ref{E:sh}) is understood in the following
sense. Any upper semi-continuous function is Borel measurable and
locally bounded from above. Hence $f|_S$ is measurable and bounded
from above on $S$. Hence the integral $\int_Sf(x)dx\in
\RR\cup\{-\infty\}$ is well defined.

(2) It is well known that for a subharmonic function
$f\not\equiv-\infty$ the integral $\int_Sf(x)dx$ is always finite.

(3) It is well known (see e.g. \cite{hormander}, Corollary 3.2.8)
that any subharmonic function $f\not\equiv-\infty$ is locally
integrable. In particular $f>-\infty$ almost everywhere. There
exists the following characterization of subharmonic functions (see
e.g. \cite{hormander}, Theorem 3.2.11). If a function $f\colon U\to
\RR\cup\{-\infty\},\, f\not\equiv -\infty,$ is subharmonic then
$\Delta f\geq 0$ in sense of generalized functions, where $\Delta
=\sum_{i=1}^n\frac{\pt^2}{\pt x_i^2}$ is the usual Laplacian on
$\RR^n$. (Notice that in order to define $\Delta f$ we use the fact
that $f$ is locally integrable.) Conversely if $U$ is a generalized
function, $\Delta U\geq 0$, then $U$ is defined by a unique
subharmonic function $u\not\equiv -\infty$. In particular twice
continuously differentiable function $f$ is subharmonic if and only
if $\Delta f\geq 0$ pointwise.
\end{remark}


Let $\Ome\subset \OO^2$ be an open subset.
\begin{definition}\label{D:psh-func0}
Let $f\colon \Ome\to \RR\cup \{-\infty\}$. The function $f$ is
called {\itshape octonionic plurisubharmonic} if

(i) $f$ is upper semi-continuous;

(ii) the restriction of $f$ to any affine octonionic line
$L\in\ca\OO\PP^1$ is subharmonic.
\end{definition}
We will denote by $P(\Ome)$ the class of all octonionic
plurisubharmonic functions in $\Ome$.

\begin{example}\label{Ex:conv-psh}
Any convex function on $\OO^2\simeq\RR^{16}$ is octonionic
subharmonic.
\end{example}

\begin{proposition}\label{psh-invar}
The class of octonionic plurisubharmonic functions is invariant
under the group $\OO^2\rtimes SL_2(\OO)$.
\end{proposition}
{\bf Proof.} This follows immediately from Lemma \ref{m-sch2}. \qed

\begin{proposition}\label{P:psh-sh}
Any octonionic plurisubharmonic function is subharmonic.
\end{proposition}
{\bf Proof.} Let $f$ be an octonionic plurisubharmonic function. We
have to show that for any $x_0$ and any sphere $S$ centered in $x_0$
(and both contained in our domain) one has
$$f(x_0)\leq \int_Sf(x)dx.$$
Without loss of generality we may assume that $x_0=0$ and $S$ has
radius 1. Then we have the equality
\begin{eqnarray}\label{E:integration}
\int_Sf(x)dx=\int_{\OO\PP^1}\left(\int_{y\in S(L)}f(y)dy\right)dL
\end{eqnarray}
where $S(L)$ denotes the unit sphere in $L$, and $dL$ is the only
$Spin(9)$-invariant probability measure on $\OO\PP^1$. The equality
(\ref{E:integration}) follows from the uniqueness of
$Spin(9)$-invariant probability measure on $S$.

Then we have
\begin{eqnarray*}
\int_Sf(x)dx=\int_{\OO\PP^1}\left(\int_{y\in
S(L)}f(y)dy\right)dL\geq \int_{\OO\PP^1}f(x_0)dL=f(x_0).
\end{eqnarray*}
\qed
\begin{corollary}\label{C:psh-integrab}
Any octonionic plurisubharmonic function $\not \equiv-\infty$ is
locally integrable.
\end{corollary}
{\bf Proof.} Since any subharmonic function $\not \equiv-\infty$ is
locally integrable (see e.g. \cite{hormander}, Corollary 3.2.8) the
result follows from Proposition \ref{P:psh-sh}. \qed

\begin{proposition}\label{P:smooth-psh}
Let $f\in C^2(\Ome)$. Then $f$ is octonionic plurisubharmonic if and
only if the matrix $\left(\dfq\right)$ is non-negative definite
pointwise.
\end{proposition}
{\bf Proof.} Let us start with the following elementary observation.
Let $b$ be a quadratic form on a Euclidean space $L$. Then
\begin{eqnarray}\label{E:integral}
Tr(b)=\dim L\cdot \int_{y\in S(L)}b(y)dy
\end{eqnarray}
where $S(L)$ denotes the unit sphere of $L$, and $dy$ is the
rotation invariant probability measure on $S(L)$.

Let us fix now $z\in \Ome$. Let us denote by $b_z$ the usual (real)
Hessian of the function $f$ at the point $z$. Let us fix also an
affine octonionic line $L$ passing through $z$. We have to show that
$\Delta_L(f|_L)|_z\geq 0$ where $\Delta_L$ denotes the Laplacian on
the line $L$. But $$\Delta_L(f|_L)|_z=2Tr(b_z|_L).$$ Using this and
(\ref{E:integral}) we get
\begin{eqnarray}\label{E:integral2}
\Delta_L(f|_L)|_z=16\int_{S(L)}b_z(y)dy.
\end{eqnarray}
We may and will assume that $L=z+\{\xi\cdot x|\,x\in \OO\}$ for some
fixed $\xi=\left[\begin{array}{c}
                    a\\
                    1
                    \end{array}\right]\in \OO^2$. Then
\begin{eqnarray*}\label{E:integral3}
\Delta_L(f|_L)|_z=\frac{16}{|\xi|^2}\int_{x\in S^7}b_z(\xi\cdot
x)dx\overset{\mbox{Lemma }\ref{L:la1}}{=}\\
\frac{1}{|\xi|^2} Re\left(\xi^*\left(\frac{\pt^2f(z)}{\pt\bar q_i\pt
q_j}\right)\xi\right).
\end{eqnarray*}
Thus we have shown that
\begin{eqnarray}\label{E:integral4}
\Delta_L(f|_L)|_z=\frac{1}{|\xi|^2}
Re\left(\xi^*\left(\frac{\pt^2f(z)}{\pt\bar q_i\pt
q_j}\right)\xi\right).
\end{eqnarray}
The last identity obviously implies the proposition. \qed



Let $dq$ denote the standard Lebesgue measure on $\OO^2$. This
choice will allow us identify functions with measures via $u\mapsto
u\cdot dq$. Thus we will not distinguish these two notions.

A {\itshape generalized function} with values in $\ch_2(\OO)$, by
definition, is a continuous linear $\RR$-valued functional on the
space $C^\infty_0(\Ome,\ch_2(\OO))$ of infinitely smooth compactly
supported functions on $\Ome$ with values in $\ch_2(\Ome)$. This
space will be denoted by $C^{-\infty}(\Ome,\ch_2(\OO))$. It is
equipped with the weak topology. We have an imbedding of the space
$L^1_{loc}(\Ome,\ch_2(\OO))$ of locally integrable
$\ch_2(\OO)$-valued functions into $C^{-\infty}(\Ome,\ch_2(\OO))$
which is given by
$$M\mapsto [\Phi\mapsto \int_\Ome Re(Tr(M\cdot \Phi))dq]$$
where $Tr$ denotes the sum of diagonal elements. (Note also that any
octonionic matrices $A$ and $B$ satisfy $Re(Tr(AB))=Re(Tr(BA))$.
\begin{definition}
Let $M\in C^{-\infty}(\Ome,\ch_2(\OO))$. We say that $M$ is
{\itshape non-negative} if for any infinitely smooth compactly
supported function $\Phi\colon \Ome\to \ch_2(\OO)$ such that
$\Phi(q)\geq 0$ for any $q\in\Ome$ one has
$$M(\Phi)\geq 0.$$
\end{definition}
\begin{remark}
If $M$ is a continuous $\ch_2(\OO)$-valued function then it is
non-negative in sense of generalized functions if and only if it is
pointwise non-negative. This follows from the observation that a
matrix $A\in \ch_2(\OO)$ is non-negative definite if and only if for
any non-negative definite matrix $B$ one has $Re(Tr(AB))\geq 0$.
\end{remark}


\begin{proposition}\label{P:sums-maximums}
(i) A linear combination of  octonionic plurisubharmonic functions
with positive coefficients is octonionic plurisubharmonic.

(ii) Maximum of two octonionic plurisubharmonic functions is
octonionic plurisubharmonic.
\end{proposition}
{\bf Proof.} This follows immediately from the corresponding
properties of subharmonic functions. \qed

Let us fix a sequence $\{\delta_n\}$ of smooth functions
approximating the $\delta$-function at the origin 0. More precisely
for any $n\in \NN$ we fix  a function $\delta_n\colon \OO^2\to \RR$
which satisfies:

(i) $\delta _n\in C^\infty(\OO^2)$;

(ii) $\delta_n\geq 0$;

(iii) $\int_{\OO^2}\delta_n(x)dx=1$;

(iv) the support $\supp(\delta_n)$ is contained in the ball of
radius $\frac{1}{n}$ centered at the origin.

\begin{proposition}\label{P:psh-approx-smooth}
Let $f\in P(\Ome)$, $f\not\equiv -\infty$ (in particular $f\in
L^1_{loc}(\Ome)$ by Corollary \ref{C:psh-integrab}).

(i) Then $f\ast\delta_n$ is infinitely smooth octonionic
plurisubharmonic function and $f\ast\delta_n\to f$ in $L^1_{loc}$.

(ii) If moreover $f\in P(\Ome)\cap C(\Ome)$ then $f\ast\delta_n\to
f$ uniformly on compact subsets of $\Ome$.
\end{proposition}
{\bf Proof.} Part (ii) in obvious. Let us prove part (i).  It is
standard (and easy to see) that for any $g\in L^1_{loc}$,
$g\ast\delta_n\to g$ weakly in sense of measures on every compact
subset of $\Ome$, in particular in the sense of generalized
functions. Observe now that in our situation the sequence $\{f\ast
\delta_n\}$ is a sequence of subharmonic functions (by Proposition
\ref{P:psh-sh}) with a uniform upper bound on every compact sunset
of $\Ome$. Since this sequence converges to $f$ in sense of
distributions, it convergence to $f$ in $L^1_{loc}$ by the general
result on subharmonic functions (\cite{hormander}, Theorem 3.2.12).
Proposition is proved. \qed

We will need the following well known fact on subharmonic functions
(see e.g. \cite{hormander}, Theorem 3.2.12).
\begin{lemma}\label{L:sh-converg}
Let $U\subset \RR^n$ be an open connected subset. Let $\{u_j\}$ be a
sequence of subharmonic functions in $U$ which have a uniform upper
bound on every compact subset of $U$.

(i) Then either $\{u_j\}\to -\infty$ uniformly on every compact
subset of $U$, or else there is a subsequence $\{u_{j_k}\}$ which
converges in $L^1_{loc}(U)$.

(ii) If $u_j\not\equiv -\infty$ for every $j$, and $\{u_j\}$
converges to $F\in C^{-\infty}(U)$ in the sense of generalized
functions, then $F$ is given by a subharmonic function $f\not\equiv
-\infty$ and $u_j\to f$ in $L^1_{loc}(U)$.
\end{lemma}
{\bf Proof.} See \cite{hormander}, Theorem 3.2.12. \qed

\begin{proposition}\label{P:sh-psh-invariance}
Let $\Ome\subset \OO^2$ be an open subset. Let $f\colon\Ome\to
\RR\cup \{-\infty\}$ be a function. Assume that for every $\phi\in
\OO^2\rtimes SL_2(\OO)$ the function $f\circ \phi$ is subharmonic in
$\phi^{-1}(\Ome)$. Then $f$ is octonionic plurisubharmonic.
\end{proposition}
{\bf Proof.} The proof is an easy modification of the proof of
Theorem 4.1.7 in \cite{hormander}. Let $z=(z_1,z_2)\in \Ome$. The
the function $f(z_1+w_1,z_2+\eps w_2)$ is subharmonic in $w$ by
hypothesis for small $\eps >0$. Hence
$$f(z)\leq \int_{|\zeta|=1} f(z_1+\zeta_1,z_2+\eps\zeta_2)d\zeta.$$
Since $f$ is upper semi-continuous and locally bounded from above,
the Fatou lemma implies as $\eps\to 0$ that
$$f(z)\leq \int_{|\zeta|=1}f(z_1+\zeta_1,z_2)d\zeta.$$
The last inequality and Theorem 3.2.3 of \cite{hormander} imply that
the function $z_1\mapsto f(z_1,z_2)$ is subharmonic. The
subharmonicity of restrictions of $f$ to other octonionic lines
follows from the hypothesis and the transitivity of the action of
$Spin(9)$ on $\OO\PP^1$. \qed
\begin{theorem}\label{T:oct-psh-convergence}
Let $\Ome\subset \OO^2$ be an open connected subset. Let $\{f_n\}$
be a sequence of octonionic plurisubharmonic functions in $\Ome$
which is uniformly bounded from above on every compact subset of
$\Ome$.

(i) Then either $f_n\to -\infty$ uniformly on every compact subset
of $\Ome$, or else there is a subsequence $\{f_{n_k}\}$ which
converges in $L^1_{loc}(\Ome)$.

(ii) If $f_n\not\equiv -\infty$ for all $n$, and $\{f_n\}$ converges
in the sense of generalized functions to $F\in C^{-\infty}(\Ome)$,
then $F$ is defined by an octonionic plurisubharmonic function
$f\not\equiv -\infty$ and $f_n\to f$ in $L^1_{loc}(\Ome)$.
\end{theorem}
{\bf Proof.} Both statements follow immediately from the
corresponding statements on subharmonic functions (Lemma
\ref{L:sh-converg}) using Propositions \ref{psh-invar},
\ref{P:psh-sh}, and \ref{P:sh-psh-invariance}. \qed

\begin{proposition}\label{P:oct-psh-charact}
Let $f\colon \Ome\to \RR\cup\{-\infty\}$ be a function such that
$f\not\equiv -\infty$. Then $f$ is octonionic plurisubharmonic if
and only if it satisfies the following conditions:

(i) $f$ is upper semi-continuous;

(ii) $f$ is locally integrable;

(iii) $\left(\dfq\right)\geq 0$ in the sense of generalized
functions.
\end{proposition}
{\bf Proof.} Actually it remains to show that if a function
$f\not\equiv -\infty$ is upper semi-continuous and locally
integrable then $f$ is octonionic plurisubharmonic if and only if
$\left(\dfq\right)\geq 0$ in the sense of generalized functions.

Let us fix a sequence $\{\delta_n\}$ approximating the
$\delta$-function at 0 as above. Then $f\ast \delta_n$ is infinitely
smooth and converges to $f$ weakly in sense of measures since $f\in
L^1_{loc}(\Ome)$. Let us assume first that $f$ is octonionic
plurisubharmonic. Then $f\ast\delta_n$ is octonionic
plurisubharmonic by Proposition \ref{P:psh-approx-smooth}(i). Hence
by Proposition \ref{P:smooth-psh} $\left(\frac{\pt^2(f\ast
\delta_n)}{\pt\bar q_i\pt  q_j}\right)\geq 0$ pointwise. But
obviously
$$\frac{\pt^2(f\ast\delta_n)}{\pt\bar q_i\pt q_j}=\frac{\pt^2f}{\pt\bar q_i\pt 
q_j}\ast\delta_n.$$ Hence $\left(\frac{\pt^2(f\ast \delta_n)}{\pt\bar
q_i\pt  q_j}\right)\to \left(\frac{\pt^2f}{\pt\bar q_i\pt 
q_j}\right)$ is the sense of generalized functions. Hence
$\left(\frac{\pt^2f}{\pt\bar q_i\pt q_j}\right)\geq 0$.

Conversely let us assume that $\left(\frac{\pt^2f}{\pt\bar q_i\pt 
q_j}\right)\geq 0$. Then $\left(\frac{\pt^2(f\ast\delta_n)}{\pt\bar
q_i\pt q_j}\right)=\left(\frac{\pt^2f}{\pt\bar q_i\pt
q_j}\right)\ast\delta_n\geq 0.$ Since $f\ast \delta_n$ is infinitely
smooth, it is octonionic plurisubharmonic by Proposition
\ref{P:smooth-psh}. Since $f\ast\delta_n\to f$ in the sense of
generalized functions, the function $f$ is octonionic
plurisubharmonic by Theorem \ref{T:oct-psh-convergence}(ii). \qed

\subsection{An analogue of the Aleksandrov and Chern-Levine-Nirenberg
theorems.}\label{Ss:cln} In this section we will denote for
brevity the octonionic Hessian $\left(\duq\right)$ by $\pt^2u$.

Let us define a 3-linear functional $\tau$ on triples of
infinitely smooth compactly supported $\RR$-valued functions on
$\OO^2$ by
$$\tau(f_0,f_1,f_2)=\int_{\OO^2}f_0\det(\pt^2 f_1,\pt^2 f_2)dq$$
where $dq$ is the standard Lebesgue measure. Later on we will use
the following technical lemma.
\begin{lemma}\label{L:psh1}
$\tau$ is symmetric with respect to all 3 arguments $f_0,f_1,f_2$.
\end{lemma}
{\bf Proof.} It is clear that $\tau$ is invariant with respect to
$f_1$ and $f_2$. It is enough to show that
\begin{eqnarray}\label{psh2}
\tau(f_0,f_1,f_2)=\tau(f_1,f_0,f_2).
\end{eqnarray}
It will be more convenient to prove (\ref{psh2}) under slightly
more general assumptions: we will assume that $f_0,f_2\in
C^\infty_0(\OO^2)$, and $f_1\in C^{-\infty}(\OO^2)$. First let us
prove (\ref{psh2}) for $f_1=\delta_{\{q_1=0\}}$, i.e. $f_1$ is the
$\delta$-function of the octonionic line $\{(0,x)|\, x\in \OO\}$.
We have
\begin{eqnarray*}
\pt^2f_1=\left[\begin{array}{cc}
                \Delta_1\delta_{\{q_1=0\}}&0\\
                0&0
                \end{array}\right],\\
\det(\pt^2f_1,\pt^2f_2)=\frac{1}{2}\Delta_2f_2\cdot
\Delta_1\delta_{\{q_1=0\}}
\end{eqnarray*}
where $\Delta_i$ denotes the usual Laplacian with respect to the
$i$-th octonionic variable.

Then
\begin{eqnarray}\label{psh3}
\tau(f_0,f_1,f_2)=\frac{1}{2} \int_{\OO^2}f_0\cdot
\Delta_2f_2\cdot\Delta_1\delta_{\{q_1=0\}}=\frac{1}{2}
\int_{\{q_1=0\}}\Delta_1(f_0\cdot \Delta_2 f_2)dq_2.
\end{eqnarray}
On the other hand
\begin{eqnarray*}
\tau(f_1,f_0,f_2)=\int_{\OO^2}f_1\cdot
\det(\pt^2f_0,\pt^2f_2)dq=\int_{\{q_1=0\}}\det(\pt^2f_0,\pt^2f_2)dq_2=\\
\int_{\{q_1=0\}}\det\left(\left[\begin{array}{cc}
                              \Delta_1f_0&\frac{\pt^2f_0}{\pt\bar
                              q_1\pt q_2}\\
                              \frac{\pt^2f_0}{\pt\bar
                              q_2\pt q_1}&\Delta_2 f_0
                              \end{array}\right],
                              \left[\begin{array}{cc}
                              \Delta_1f_2&\frac{\pt^2f_2}{\pt\bar
                              q_1\pt q_2}\\
                              \frac{\pt^2f_2}{\pt\bar
                              q_2\pt q_1}&\Delta_2 f_2
                              \end{array}\right]
                              \right)dq_2=\\
\frac{1}{2}\int_{\{q_1=0\}}\left(\Delta_1f_0\cdot
\Delta_2f_2+\Delta_1f_2\cdot \Delta_2
f_0-2Re\left(\frac{\pt^2f_0}{\pt\bar q_1\pt
q_2}\cdot\frac{\pt^2f_2}{\pt\bar q_2\pt
q_1}\right)\right)dq_2.\end{eqnarray*} After integration by parts
in the second summand we obtain
\begin{eqnarray}\label{psh4}
\tau(f_1,f_0,f_2)=\frac{1}{2}\int_{\{q_1=0\}}\left(\Delta_1f_0\cdot
\Delta_2f_2+f_0\cdot \Delta_1\Delta_2f_2
-2Re\left(\frac{\pt^2f_2}{\pt\bar q_2\pt q_1}\cdot
\frac{\pt^2f_0}{\pt\bar q_1\pt q_2} \right)\right)dq_2
\end{eqnarray}
Let us integrate the third summand in the last expression:
\begin{eqnarray*}
\int_{\{q_1=0\}}Re\left(\frac{\pt^2f_2}{\pt\bar q_2\pt q_1}\cdot
\frac{\pt^2f_0}{\pt\bar q_1\pt q_2}
\right)dq_2=\\
\int_{\{q_1=0\}}Re\left(\left(\sum_{t=0}^7e_t\frac{\pt^2f_2}{\pt
x^{t}_2\pt  q_1}\right)\cdot
\left(\sum_{s=0}^7\frac{\pt^2f_0}{\pt x^{s}_2\pt\bar
q_1}\bar e_s\right)\right)dq_2\overset{\mbox{by parts}}{=}\\
-\sum_{s,t=0}^7\int_{\{q_1=0\}}Re\left(\left(e_t\frac{\pt^2}{\pt
x^{s}_2\pt x^{t}_2}\left(\frac{\pt f_2}{\pt q_1}\right)\right)\cdot \left(\frac{\pt f_0}{\pt\bar
q_1}\bar e_s\right)\right)dq_2\overset{\mbox{Lemma }\ref{L:oalg1}(i)}{=}\\
-\sum_{s,t=0}^7\int_{\{q_1=0\}}Re\left(\left(\bar e_s\left(e_t\frac{\pt^2}{\pt
x^{s}_2\pt x^{t}_2}\left(\frac{\pt f_2}{\pt q_1}\right)\right)\right)\cdot\frac{\pt f_0}{\pt\bar
q_1}\right)dq_2\overset{\mbox{Lemma }\ref{L:oalg1}(iii)}{=}\\
-\sum_{s=0}^7\int_{\{q_1=0\}}Re\left(\left(\frac{\pt^2}{ (\pt
x^{s}_2)^2}\left(\frac{\pt f_2}{\pt q_1}\right)\right)\cdot
\frac{\pt f_0}{\pt\bar q_1}\right)dq_2-\\
\sum_{0\leq s<t\leq
7}\int_{\{q_1=0\}}Re\left(\left(\left(\frac{\pt^2}{\pt x^{s}_2\pt
x^{t}_2}\left(\frac{\pt f_2}{\pt  q_1}\right)\right) (\bar e_t
e_s+\bar e_s e_t)\right)\cdot\frac{\pt f_0}{\pt\bar
q_1}\right)dq_2=\\
-\int_{\{q_1=0\}}Re\left(\frac{\pt(\Delta_2 f_2)}{\pt
q_1}\cdot\frac{\pt f_0}{\pt\bar q_1}\right)dq_2.
\end{eqnarray*}
Substituting the last expression back to (\ref{psh4}) we obtain
\begin{eqnarray*}
\tau(f_1,f_0,f_2)=\frac{1}{2}\int_{\{q_1=0\}}\left(\Delta_1f_0\cdot
\Delta_2f_2+f_0\cdot
\Delta_1\Delta_2f_2+2Re\left(\frac{\pt(\Delta_2 f_2)}{\pt
q_1}\cdot\frac{\pt f_0}{\pt\bar q_1}\right)\right)dq_2=\\
\frac{1}{2}\int_{\{q_1=0\}}\frac{\pt^2}{\pt\bar q_1\pt
q_1}\left(f_0\cdot\Delta_2f_2\right)dq_2=
\frac{1}{2}\int_{\{q_1=0\}}\Delta_1\left(f_0\cdot\Delta_2f_2\right)dq_2\overset{\mbox{by
} (\ref{psh3})}{=} \tau(f_0,f_1,f_2).
\end{eqnarray*}
Thus the equality (\ref{psh3}) is proved for
$f_1=\delta_{\{q_1=0\}}$. Next using Proposition \ref{P:ohs1} and
the fact that the group $\OO^2\rtimes Spin(9)$ acts transitively on
$\ca\OO\PP^1$ we conclude that the equality (\ref{psh2}) for
$f_1=\delta_L$ for any $L\in\ca\OO\PP^1$. But by Corollary
\ref{C:radon3} linear combinations of such $\delta$-functions are
dense in $C^{-\infty}(\OO^2)$, hence the equality (\ref{psh2}) is
proved for any generalized function $f_1$. \qed


For any matrix-valued function $F\colon \Omega\to \ch_2(\OO)$ let
us denote by $||F||_{L^1(\Omega)}$ its $L^1$-norm, i.e. the sum of
the $L^1$-norms of the elements of this matrix. It is easy to see
that if $F$ takes values in non-negative definite matrices, then
\begin{eqnarray}\label{E:cln0.1}
||F||_{L^1(\Omega)}\leq
2||\det(F,I_2)||_{L^1(\Omega)}=2\int_\Omega\det(F,I_2)
\end{eqnarray}
where $I_2\in \ch_2(\OO)$ is the identity matrix.
\begin{lemma}\label{L:cln1}
For any compact subset $K\subset \Omega$ and any its compact
neighborhood $K'\subset\Omega$ there exists a constant $C$ such
that for any $f\in P(\Omega)\cap C^2(\Omega)$ one has
\begin{eqnarray}\label{E:cln2}
||(\pt^2f)||_{L^1(K)}\leq C||f||_{L^\infty(K')};\\\label{E:cln3}
||\det(\pt^2f)||_{L^1(K)}\leq C||f||^2_{L^\infty(K')}.
\end{eqnarray}
\end{lemma}
{\bf Proof.} Let us fix a non-negative function $\gamma\in
C^\infty_0(\Omega)$ which is equal to 1 in a neighborhood of $K$
and vanishes on a neighborhood of the closure
$\overline{\Omega\backslash K'}$. Then using (\ref{E:cln0.1}) one
has
\begin{eqnarray*}
||(\pt^2f)||_{L^1(K)}\leq 2\int_\Omega\gamma
\det(\pt^2f,I_2)=\\
\int_\Omega \gamma \cdot \Delta f=\int_\Omega f\Delta \gamma \leq
C||f||_{L^\infty(K')}.
\end{eqnarray*}
Thus the inequality (\ref{E:cln2}) is proved. Let us prove
(\ref{E:cln3}). We have
\begin{eqnarray*}
||\det(\pt^2f)||_{L^1(K)}\leq \int_\Omega\gamma
\det(\pt^2f)\overset{\mbox{Lemma }\ref{L:psh1}}{=} \int_\Omega
f\det(\pt^2\gamma,\pt^2f)\leq\\C' ||f||_{L^\infty(K')}\cdot
||(\pt^2f)||_{L^1(\supp \gamma)}\leq C''||f||_{L^\infty(K')}^2
\end{eqnarray*}
where the last inequality follows from (\ref{E:cln2}). Lemma is
proved. \qed

\begin{corollary}\label{C:cln4}
For any compact subset $K\subset \Omega$ and any its compact
neighborhood $K'\subset \Omega$ there exists a constant $C$ such
that for any $f,g\in P(\Omega)\cap C^2(\Omega)$ and any $\phi\in
C^2(\Omega)$ with $\supp\phi\subset K$ one has
\begin{eqnarray}\label{E:cln5}
\big|\int_\Omega \phi\left(\det(\pt^2f)-\det(\pt^2g)\right)\big|
\leq C||f-g||_{L^\infty(K)}
\left(||f||_{L^\infty(K')}+||g||_{L^\infty(K')}\right)\cdot
||\phi||_{C^2(\Omega)}.
\end{eqnarray}
\end{corollary}
{\bf Proof.} We have
\begin{eqnarray*}
\big|\int_\Omega
\phi\left(\det(\pt^2f)-\det(\pt^2g)\right)\big|\leq\\
\big|\int\phi\det(\pt^2(f-g),\pt^2f)\big|+\big|\int\phi\det(\pt^2g,\pt^2(f-g))\big|
\overset{\mbox{Lemma } \ref{L:psh1}}{=}\\
\big|\int(f-g)\det(\pt^2f,\pt^2\phi)\big|+\big|\int(f-g)\det(\pt^2g,\pt^2\phi)\big|\leq\\
C||f-g||_{L^\infty(K)}||\phi||_{C^2(\Omega)}\left(||(\pt^2f)||_{L^1(K)}+||(\pt^2g)||_{L^1(K)}\right)
\overset{\mbox{Lemma }\ref{L:cln1}}{\leq}\\
C'||f-g||_{L^\infty(K)}||\phi||_{C^2(\Omega)}\left(||f||_{L^\infty(K')}+||g||_{L^\infty(K')}\right).
\end{eqnarray*}
Corollary is proved. \qed

\begin{proposition}
Let $\Ome\subset \OO^2$ be an open subset. Assume that a sequence
$\{f_n\}\subset P(\Ome)$, $f_n\not\equiv -\infty$ for any $n$,
converges in $L^1_{loc}(\Ome)$ to an octonionic plurisubharmonic
function $f\not\equiv -\infty$. Then one has a week convergence of
$\ch_2(\OO)$-valued measures
$$\left(\frac{\pt^2f_n}{\pt\bar q_i \pt q_j}\right)\to
\left(\dfq\right).$$
\end{proposition}
\def\dun{\frac{\pt^2f_n}{\pt\bar q_i\pt q_j}}
{\bf Proof.} By Lemma \ref{L:cln1} the measures $\{(\dun)\}$ are
uniformly locally bounded in $\Omega$. Hence choosing a subsequence
if necessary we may assume that  this sequence of measures converges
weakly to an $\ch_2(\OO)$-valued measure $(\nu _{\bar i j})$. We have
to prove that $\nu _{\bar i j}= \dfq$. To see it let us fix an
arbitrary function $\phi\in C^{\infty}_0(\Omega)$. Then
$$\int _{\Omega} \dun \cdot \phi dq =
 \int _{\Omega} f_n \cdot \frac{\partial ^2 \phi}{\partial\bar q_i \partial q _j}
 dq
\to \int _{\Omega} f \cdot \frac{\partial ^2 \phi}{\partial\bar q_i
\partial q _j} dq= \int _{\Omega}\frac{\partial ^2 f}{\partial\bar q_i
\partial q _j}\phi \cdot dq ,$$ where the first and the last equalities
are obtained by integration by parts. The result follows. \qed

\begin{proposition}\label{P:cln0}
Let $f\in P(\Omega)\cap C(\Omega)$. There exists a unique measure
on $\Ome$ denoted by $\det(\pt^2f)$ satisfying the following
property: for any sequence $\{f_n\}\subset P(\Omega)\cap
C^2(\Omega)$ such that $f_n\to f$ uniformly on compact subsets of
$\Ome$ one has $\det(\pt^2f_n)\to \det(\pt^2f)$ weakly in sense of
measures. This measure $\det(\pt^2 f)$ is non-negative and has the
obvious interpretation when $f\in C^2(\Omega)$.
\end{proposition}
{\bf Proof.} By Lemma \ref{L:cln1} for any compact subset
$K\subset \Omega$ and any its compact neighborhood $K'$ there
exists a constant $C$ such that for any $n$
$$||\det(\pt^2f_n)||_{L^1(K)}\leq C||f_n||^2_{L^\infty(K')}.$$
Hence the sequence of measures $\det(\pt^2f_n)|_K$ on $K$ has a
compact closure in the weak topology. It remains to show that this
sequence has at most one limiting point. By Corollary \ref{C:cln4}
we have for any $\phi\in C^\infty_0(\Omega)$ with
$\supp\phi\subset K$
\begin{eqnarray}
\big|\int_\Ome\phi\left(\det(\pt^2f_n)-\det(\pt^2f_m)\right)\big|\leq
C||\phi||_{C^2(\Omega)}||f_n-f_m||_{L^\infty(K)}(||f_n||_{L^\infty(K')}+||f_m||_{L^\infty(K')}).
\end{eqnarray}
Obviously the last expression tends to 0 as $m,n\to\infty$. This
implies also that the limiting measure is unique and it is
independent of the approximating sequence $\{f_n\}$. Let us denote
it temporarily by $\mu$.

Obviously the limiting measure $\mu$ is non-negative. It remains to
check that if  $f\in C^2(\Omega)$ then $\mu=\det(\pt^2f)$. We have
shown that $\mu$ is independent of an approximating sequence. Taking
the constant sequence equal to $f$ we conclude the result. \qed

The next lemma generalizes Lemma \ref{L:cln1} and Corollary
\ref{C:cln4} to the class of all continuous plurisubharmonic
functions.
\begin{lemma}\label{L:cln10}
For any compact subset $K\subset \Omega$ and any its compact
neighborhood $K'\subset\Omega$ there exists a constant $C$ such
that

(i) for any $f\in P(\Omega)\cap C(\Omega)$ one has
\begin{eqnarray}\label{E:cln20}
||(\pt^2f)||_{L^1(K)}\leq C||f||_{L^\infty(K')},\\\label{E:cln30}
||\det(\pt^2f)||_{L^1(K)}\leq C||f||^2_{L^\infty(K')};
\end{eqnarray}

(ii) for any $f,g\in P(\Omega)\cap C(\Omega)$ and any $\phi\in
C^2(\Omega)$ with $\supp\phi\subset K$ one has
\begin{eqnarray}\label{E:cln50}
\big|\int_\Omega \phi\left(\det(\pt^2f)-\det(\pt^2g)\right)\big|
\leq C||f-g||_{L^\infty(K)}
\left(||f||_{L^\infty(K')}+||g||_{L^\infty(K')}\right)\cdot
||\phi||_{C^2(\Omega)}.
\end{eqnarray}
\end{lemma}
{\bf Proof.} This lemma follows immediately from Lemma \ref{L:cln1}
and Corollary \ref{C:cln4} using Proposition \ref{P:cln0} and the
approximation of $f$ by smooth plurisubharmonic functions as in
Proposition \ref{P:psh-approx-smooth}. \qed

\begin{theorem}\label{T:cln-main}
Let $\Omega\subset\OO^2$ be an open subset.  Let a sequence
$\{f_n\}\subset P(\Omega)\cap C(\Ome)$ converges uniformly on
compact subsets of $\Ome$ to a function $f$. Then $f\in
P(\Omega)\cap C(\Omega)$ and $\det(\pt^2f_n)\to \det(\pt^2f)$ weakly
in sense of measures.
\end{theorem}
{\bf Proof.} This theorem follows from Lemma \ref{L:cln10} exactly
in the same way as Proposition \ref{P:cln0} followed from Lemma
\ref{L:cln1} and Corollary \ref{C:cln4}. \qed

From Proposition \ref{P:cln0} and Theorem \ref{T:cln-main} one can
easily deduce the following 'mixed' version of these results
generalizing both of them. This version uses the notion of mixed
determinant introduced in Section \ref{S:la}.
\begin{theorem}\label{T:cln-mixed}
Let $\Ome\subset \OO^2$ be an open subset. For any $f,g\in
P(\Ome)\cap C(\Ome)$ there exists a non-negative measure denoted
by $\det(\pt^2f,\pt^2g)$. It satisfies and is uniquely
characterized by the following two properties:

(i) if $f,g\in C^2(\Ome)$ then the meaning is obvious;

(ii) if one has two sequences $\{f_n\},\{g_n\}\subset P(\Ome)\cap
C(\Ome)$ such that $f_n\to f,\, g_n\to g$ uniformly on compact
subsets of $\Ome$ then $\det\left(\pt^2 f_n ,\pt^2 g_n\right)\to
\det\left(\pt^2 f ,\pt^2 g\right)$ weakly in the sense of measures.
\end{theorem}
Note that non-negativity of measures follows from Lemma
\ref{L:mix-det-posit}.

\subsection{A B{\l}ocki type theorem.}\label{Ss:blocki}
\setcounter{theorem}{0}
\begin{theorem}\label{T:blocki1}
For any $u,v\in P(\Ome)\cap C(\Ome)$ one has
\begin{eqnarray}\label{blocki-for}
\det(\pt^2(\max\{u,v\}))=\det(\pt^2(\max\{u,v\}),\pt^2u+\pt^2v)-\det(\pt^2u,\pt^2v).
\end{eqnarray}
\end{theorem}
{\bf Proof.} The argument is very close to the original B{\l}ocki's
argument \cite{blocki}. By continuity of both sides in
(\ref{blocki-for}) we may assume that $u,v$ are smooth. Let
$\chi:\RR \to [0,\infty)$ be a smooth function such that $\chi(x)=0$
if $x\leq -1$, $\chi(x)=x$ if $x\geq 1$, and $0\leq \chi ' \leq 1,\,
\chi '' \geq 0$ everywhere. Define
$$\psi_j:=v+ \frac{1}{j}\chi(j(u-v)),$$
$$\alpha:=u-v,$$
$$w:=\max\{u,v\}.$$ It is easy to see that $\psi_j\downarrow w$
uniformly on compact subsets and monotonically as $j\to \infty$.
\begin{lemma}\label{308}
$$\left(\frac{\chi(j\alpha)}{j}\right)_{\bar p q}= \chi
'(j\alpha)\cdot \alpha_{\bar p q}+ j\chi ''(j\alpha)\alpha_{\bar
p}\alpha_q.$$
\end{lemma}
{\bf Proof.} Since $\alp$ is a real valued function we have
$$\left(\frac{\chi(j\alpha)}{j}\right)_{p\bar q}=\frac{1}{j}
\sum_{l,m=0}^7 e_l(\chi (j\alpha))_{x_q^l x_p^m}\bar e_m=$$
$$\sum_{l,m=0}^7 e_l\left( \chi '(j\alpha)\cdot \alpha_{x_p^l}\right)
_{x_q^m}\bar e_m = \sum_{l,m=0}^7 e_l \left(j \chi ''(j\alpha)
\cdot \alpha _{x_p^l}\alpha _{x_q^m} + \chi '(j\alpha )
\alpha_{x_p^l x_q^m}\right) \bar e_m=$$
$$\chi '(j\alpha) \alpha_{\bar p q} + j \chi ''(j\alpha)\alpha
_{\bar p}\alpha _q.$$ \qed

Thus from Lemma \ref{308} we obtain
$$ (\psi_j)_{\bar p q}=v_{\bar p q}+ \chi' (j\alpha)(u-v)_{\bar p q}+
j \chi ''(j\alpha)\alpha _{\bar p}\alpha _q=$$
$$\chi '(j\alpha)u_{\bar p q}+ (1-\chi'(j\alpha ))v_{\bar p q}+j \chi ''(j\alpha)\alpha
_{\bar p}\alpha _q.$$ The matrix $\left(\alp_{\bar p}\alp_q\right)$
is non-negative definite. Then, since $0\leq \chi'\leq 1$ and
$\chi''\geq 0$, it follows that $\psi_j$ is plurisubharmonic. From
the definition of $\psi_j$ we have
\begin{equation}\label{2}\det(
\pt^2\psi_j)=\det(\pt^2 v)+
2\det\left(\pt^2v,\pt^2\left(\frac{\chi(j\alpha)}{j}\right)
\right)+\det\left(\pt^2\left( \frac{\chi(j\alpha)}{j}\right)
\right).
\end{equation}
We have weak convergence
\begin{eqnarray}
& & \det(\pt^2\psi_j)\to \det(\pt^2 w),\label{3}\\
& & \det\left(\pt^2v,
\pt^2\left(\frac{\chi(j\alpha)}{j}\right)\right)\to
\det\left(\pt^2(w-v), \pt^2 v\right) \label{4}.
\end{eqnarray}
Let us study the last term in (\ref{2}), namely
$\det\left(\pt^2\left( \frac{\chi(j\alpha)}{j}\right) \right)$. From
Lemma \ref{308} one gets
\begin{eqnarray}\label{E:vst1}\label{star0}
\det\left(\pt^2\left(
\frac{\chi(j\alpha)}{j}\right)\right)=\\\label{vyr} \chi
'(j\alpha)^2 \det\left(\alpha_{\bar p q}\right)+ 2j\chi
'(j\alpha)\chi''(j\alpha)\det\left((\alpha_{\bar p}\alpha_q),
(\alpha_{\bar p q})\right)+(j\chi''(j\alp))^2\det(\alp_{\bar
p}\alp_q)=\\\label{star} \chi '(j\alpha)^2 \det\left(\alpha_{\bar p
q}\right)+ 2j\chi '(j\alpha)\chi''(j\alpha)\det\left((\alpha_{\bar
p}\alpha_q), (\alpha_{\bar p q})\right)
\end{eqnarray}
since the last summand in (\ref{vyr}) vanishes. Let $\gamma:\RR \to
\RR$ be such that $\gamma' =(\chi')^2$. Then we have
\begin{lemma}\label{3010}
$$\det\left(\pt^2\left( \frac{\chi(j\alpha)}{j}\right) \right)=
\det\left(\pt^2\left(\frac{\gamma(j\alpha)}{j}\right),
\pt^2\alpha\right).$$
\end{lemma}
Let us postpone the proof of Lemma \ref{3010}, and finish the proof
of Theorem \ref{T:blocki1}. One can choose $\gamma$ so that $\gamma
(-1)=0$. Then $\frac{\gamma(jx)}{j}\downarrow \max\{0,x\}$ uniformly
on compact subsets and monotonically as $j\to \infty$. Hence
$$\det\left(\pt^2\left( \frac{\chi(j\alpha)}{j}\right) \right)\to
\det(\pt^2(w-v), \pt^2\alpha) \mbox{ weakly }.$$ This and
(\ref{2}), (\ref{3}), (\ref{4}) imply
\begin{eqnarray*}\det(\pt^2w)&=&\det(\pt^2v)
+2\det(\pt^2(w-v), \pt^2v)+\det(\pt^2(w-v), \pt^2(u-v))\\
&=& \det(\pt^2w, (\pt^2u+\pt^2v))-\det(\pt^2u, \pt^2v).
\end{eqnarray*}
This implies Theorem \ref{T:blocki1}. It remains to prove Lemma
\ref{3010}.

{\bf Proof of Lemma \ref{3010}.} We have
\begin{eqnarray*}
\left(\frac{\gamma(j\alpha)}{j}\right)_{\bar p q}&=&
\gamma'(j\alpha)\alpha_{\bar p q}+\alp_{\bar p} \cdot
(\gamma'(j\alp))_q\\
&= & (\chi'(j\alp))^2\alpha_{\bar p q}+ 2j \chi'(j\alp)
\chi''(j\alp)\alp_{\bar p}\cdot \alp_q.
\end{eqnarray*}
This and the equality (\ref{star})=(\ref{star0}) imply Lemma
\ref{3010}. \qed

\begin{corollary}\label{C:blocki3}
For any $u,v\in P(\Ome)\cap C(\Ome)$ such that $\min\{u,v\}\in
P(\Ome)$ one has
$$\det(\pt^2(\min\{u,v\}))=\det(\pt^2u)+\det(\pt^2v)-\det(\pt^2(\max\{u,v\})).$$
\end{corollary}
{\bf Proof.} Observe that $\min\{u,v\}=u+v-\max\{u,v\}$. Denote
for brevity $U:=\pt^2u,V:=\pt^2v, W:=\pt^2(\max\{u,v\})$. Then we
get
\begin{eqnarray*}
\det(\pt^2(\min\{u,v\}))=\det(U+V-W)=\\
(\det U+\det V-\det W)+\\2\det W
+2\det(U,V)-2\det(U,W)-2\det(V,W)\overset{\mbox{Theorem
}\ref{T:blocki1}}{=}\\\det U+\det V-\det W.
\end{eqnarray*}
Corollary is proved. \qed

\section{Valuations on convex subsets of
$\OO^2$.}\label{S:valuat}\setcounter{subsection}{0}\setcounter{theorem}{0}

\subsection{"Obvious" examples of valuations on
$\OO^2$.}\label{Ss:obvious-examples}
\def\vspn{Val^{Spin(9)}(\OO^2)}
We denote by $Val(\OO^2)$ the space of translation invariant
continuous valuations on convex compact subsets of $\OO^2$. We
denote by $\vspn$ the subspace of $Spin(9)$-invariant valuations.
Since the group $Spin(9)$ acts transitively on the unit sphere
$S^{15}$, the space $\vspn$ is finite dimensional by
\cite{alesker-adv-00}, Theorem 8.1. Since $-Id\in Spin(9)$ all
valuations in $\vspn$ are even, i.e. they take the same value on $K$
and $-K$ for any convex compact set $K$.

For the moment we can neither classify valuations in $\vspn$ nor
even compute the dimension of this space. The goal of this section
is to present few examples of such valuations which are natural from
the point of view of convexity and integral geometry.

For a non-negative integer $i\leq 16$ let us denote by $V_i$ the
intrinsic volume which is a continuous valuation invariant under all
isometries of the Euclidean space $\RR^{16}$. By definition, for any
$K\in\ck(\OO^2)$ the intrinsic volume $V_i(K)$ is equal (up to a
normalizing constant) to the mixed volume $V(\underset{i \mbox{
times}} {\underbrace{K, \dots,K}};\underset{16-i \mbox{ times}}
{\underbrace{D, \dots,D}})$ where $D$ is the unit Euclidean ball. We
refer to the book \cite{schneider-book} for the details on the mixed
and intrinsic volumes. Thus in particular $V_i\in \vspn$. Recall
that $V_{16}$ is proportional to the Lebesgue measure, and $V_0$ is
the Euler characteristic (which is equal to 1 on any convex compact
set).

To construct more examples let us fix $i=0,1,\dots, 8$. Define
\begin{eqnarray}\label{Ex:T-val}
T_i(K):=\int_{E\in \OO\PP^1}V_i(pr_E(K))dE
\end{eqnarray}
where $pr_E\colon \OO^2\to E$ is the orthogonal projection, and $dE$
is the probability $Spin(9)$-invariant Haar measure on $\OO\PP^1$.
It is easy to see that $T_i\in \vspn$. Moreover $T_i$ is
$i$-homogeneous: a valuation $\phi$ is called $i$-homogeneous if
$\phi(\lam K)=\lam^i\phi(K)$ for any $\lam >0$ and any set $K$.

Next fix $8\leq j\leq 16$. Define
\begin{eqnarray}\label{Ex:U-val}
U_j(K):=\int_{F\in\ca\OO\PP^1}V_{j-8}(K\cap F)dF
\end{eqnarray} where $dF$ is $\OO^2\rtimes Spin (9)$-invariant Haar
measure on $\ca\OO\PP^1$ (we do not specify the normalization of
this measure since it is irrelevant for the moment). Then $U_j\in
\vspn$. $U_j$ is $j$-homogeneous.

\begin{remark}\label{R:obv-ex1}
The valuations $T_i$ and $U_{16-i}$ are Fourier transform of each
other (up to a constant). For the notion of the Fourier transform on
valuations we refer to \cite{alesker-jdg-03} for the even case
(where this transform is called the duality operator) and to
\cite{alesker-fourier} for the general case. Notice that the
intrinsic volumes $V_i$ and $V_{16-i}$ are also Fourier transforms
of each other (up to a constant).
\end{remark}

\begin{remark}\label{R:obv-ex2}
$\vspn$ has a natural product \cite{alesker-gafa-04} making it a
commutative associative graded algebra with unit (where the unit is
the Euler characteristic). Thus taking polynomials in the above
examples of valuations we can produce more examples of
$Spin(9)$-invariant valuations. Moreover one can take convolutions
in sense of \cite{bernig-fu} of the above examples. However at
present relations of these examples to the previous ones are not
known (e.g. which of them are linearly independent? do they span
$\vspn$?).
\end{remark}

\begin{remark}\label{R:obv-ex3}
There is still yet another general construction of valuations which
is based on integration of $Spin(9)$-invariant differential forms on
the spherical cotangent bundle of $\OO^2$ with respect to the normal
cycle of a convex set. It was shown in \cite{alesker-mflds1},
Theorem 5.2.1, that this construction gives all valuations from
$\vspn$. This construction of valuations uses results of J. Fu
\cite{fu-90}, \cite{fu-94}; more details in the context of
valuations are given in \cite{alesker-fu-mflds3}. Relations of this
construction to other ones discussed in this article have not been
studied.
\end{remark}

\subsection{New examples of valuations on
$\OO^2$.}\label{Ss:new-examples} Recall that the space $\OO^2$ is
equipped with the standard Euclidean product
$$<(q_1,q_2),(z_1,z_2)>=Re(q_1\bar z_1+q_2\bar z_2).$$
For a convex compact set $K\in \ck(\OO^2)$ one defines its {\itshape
supporting functional} $h_K\colon \OO^2\to \RR$ by
$$h_K(x):=\sup_{y\in K} <x,y>.$$
Then $h_K$ is a convex 1-homogeneous function. In particular it is
octonionic plurisubharmonic (Example \ref{Ex:conv-psh}).

\begin{theorem}\label{T:valuat-main}
Fix a continuous compactly supported function $\psi$ on $\OO^2$.
Then
$$K\mapsto \int_{\OO^2}\det\left(\frac{\pt^2h_K}{\pt\bar q_i\pt q_j}\right)\cdot
\psi dq$$ is a translation invariant continuous 2-homogeneous
valuation on $\ck(\OO^2)$.
\end{theorem}
{\bf Proof.} Translation invariance is obvious. Continuity follows
from Theorem \ref{T:cln-main}. To prove the valuation property let
us observe first that if $K=K_1\cup K_2$ with $K_1,K_2,K\in
\ck(\OO^2)$ then
$$h_K=\max\{h_{K_1},h_{K_2}\},\, h_{K_1\cap
K_2}=\min\{h_{K_1},h_{K_2}\}.$$ Hence the result follows from
Theorem \ref{T:blocki1}. \qed

From Theorem \ref{T:valuat-main} and Proposition \ref{P:ohs1} we
immediately deduce the following corollary.
\begin{corollary}\label{C:val-invar}
The correspondence
$$K\mapsto \int_D\det\left(\frac{\pt^2h_K}{\pt\bar q_i\pt
q_j}\right)dq$$ where $D$ is the unit ball in $\OO^2$, is a
$Spin(9)$-invariant translation invariant continuous 2-homogeneous
valuation on $\ck(\OO^2)$.
\end{corollary}

\begin{remark}
It is not hard to see that the valuation from Corollary
\ref{C:val-invar} is {\itshape not} invariant under the group
$SO(16)$. In particular it is not proportional to the second
intrinsic volume $V_2$ on $\RR^{16}$.
\end{remark}

\end{document}